\documentclass[11pt]{article}%
\usepackage{amssymb}
\usepackage{graphicx}
\usepackage{amsthm,amsmath,amsfonts}
\usepackage{amsmath}
\usepackage{amsfonts}%
\newtheorem{theorem}{Theorem}
\newtheorem{lemma}{Lemma}
\newtheorem{proposition}{Proposition}
\newtheorem{corollary}{Corollary}
\newcounter{remarknumber}
\newenvironment{remark}{\stepcounter{remarknumber} \noindent
{\bfseries{Remark \theremarknumber.}}}{\medskip}
\begin{document}

\date{}
\title{Randomly Weighted Self-normalized L\'{e}vy Processes}
\author{P\'eter Kevei\thanks{Supported by the TAMOP--4.2.1/B--09/1/KONV--2010--0005
project.}\\Analysis and Stochastics Research Group of the Hungarian Academy of Sciences \\Bolyai Institute, Aradi v\'ertan\'uk tere 1, 6720 Szeged, Hungary and \\CIMAT, Callej\'on Jalisco S/N, Mineral de Valenciana, Guanajuato 36240, Mexico \\e-mail: \texttt{kevei@math.u-szeged.hu} \smallskip\\David M. Mason\thanks{Research partially supported by NSF Grant DMS--0503908.}\\Statistics Program, University of Delaware \\213 Townsend Hall, Newark, DE 19716, USA \\e-mail: \texttt{davidm@udel.edu}}
\maketitle

\begin{abstract}
\noindent Let $(U_{t},V_{t})$ be a bivariate L\'{e}vy process, where $V_{t}$
is a subordinator and $U_{t}$ is a L\'{e}vy process formed by randomly
weighting each jump of $V_{t}$ by an independent random variable $X_{t}$
having cdf $F$. We investigate the asymptotic distribution of the
self-normalized L\'{e}vy process $U_{t}/V_{t}$ at 0 and at $\infty$. We show
that all subsequential limits of this ratio at 0 ($\infty$) are continuous for
any nondegenerate $F$ with finite expectation if and only if $V_{t}$ belongs
to the centered Feller class at 0 ($\infty$). We also characterize when
$U_{t}/V_{t}$ has a non-degenerate limit distribution at 0 and $\infty$.

\textit{AMS Subject Classification:} MSC 60G51; MSC 60F05.

\textit{Keywords:} L\'evy process; Feller class; self-normalization; stable distributions.

\end{abstract}

\section{Introduction and statements of two main results}

We begin by defining the bivariate L\'{e}vy process $\left(  U_{t}%
,V_{t}\right)  ,$ $t\geq0,$ that will be the object of our study. Let $F$ be a
cumulative distribution function [cdf] satisfying
\begin{equation}
\int_{-\infty}^{\infty}\left\vert x\right\vert F\left(  \mathrm{d}x\right)
<\infty\label{E}%
\end{equation}
and $\Lambda$ be a L\'{e}vy measure on $\mathbb{R}$ with support in $\left(
0,\infty\right)  $ such that
\begin{equation}
\int_{0}^{1}y\Lambda(\mathrm{d}y)<\infty.\label{VV}%
\end{equation}
We define the\textit{\ L\'{e}vy function} $\overline{\Lambda}\left(  x\right)
=\Lambda\left(  x,\infty\right)  $ for $x\geq0$. Using Corollary 15.8 on page
291 of Kallenberg \cite{kallenberg} and assumptions (\ref{E}) and (\ref{VV}),
we can define via $F$ and $\Lambda$ the bivariate L\'{e}vy process $\left(
U_{t},V_{t}\right)  ,$ $t\geq0$, having the joint characteristic function
\begin{equation}
E\exp\left(  \mathrm{i}\theta_{1}U_{t}+\mathrm{i}\theta_{2}V_{t}\right)
=:\phi\left(  t,\theta_{1},\theta_{2}\right)  =\exp\left(  t\int_{(0,\infty
)}\int_{-\infty}^{\infty}\left(  e^{\mathrm{i}(\theta_{1}u+\theta_{2}%
v)}-1\right)  \Pi(\mathrm{d}u,\mathrm{d}v)\right)  ,\label{holds}%
\end{equation}
with%
\begin{equation}
\Pi\left(  \mathrm{d}u,\mathrm{d}v\right)  =F\left(  \mathrm{d}u/v\right)
\Lambda\left(  \mathrm{d}v\right)  .\label{pi}%
\end{equation}
From the form of $\phi\left(  t,\theta_{1},\theta_{2}\right)  $ it is clear
that $V_{t}$ is a driftless subordinator.\smallskip\ 

Throughout this paper $\left(  U_{t},V_{t}\right)  $, $t\geq0,$ denotes a
L\'{e}vy process satisfying (\ref{E}) and (\ref{VV}) and having joint
characteristic function (\ref{holds}).\smallskip

Now let $\left\{  X_{s}\right\}  _{s\geq0}$ be a class of i.i.d. $F$ random
variables independent of the $V_{t}$ process. We shall soon see that for each
$t\geq0$ the bivariate process%
\begin{equation}
\left(  U_{t},V_{t}\right)  \overset{\mathrm{D}}{=}\left(  \sum_{0\leq s\leq
t}X_{s}\Delta V_{s},\sum_{0\leq s\leq t}\Delta V_{s}\right)  ,\label{WL}%
\end{equation}
where $\Delta V_{s}=V_{s}-V_{s-}$. Notice that in the representation
(\ref{WL}) each jump of $V_{t}$ \ is weighted by an independent $X_{t}$ so
that $U_{t}$ can be viewed as a randomly weighted L\'{e}vy process.\smallskip

Here is a graphic way to picture this bivariate process. Consider $\Delta
V_{s}$ as the intensity of a random shock to a system at time $s>0$ and
$X_{s}\Delta V_{s}$ as the cost of repairing the damage that it causes. Then
$V_{t}$, $U_{t}$ and $U_{t}/V_{t}$ represent, respectively, up to time $t$,
the total intensity of the shocks, the total cost of repair and the average
cost of repair with respect to shock intensity. For instance, $\Delta V_{s}$
can represent a measure of the intensity of a tornado that comes down in a
Midwestern American state at time $s$ during tornado season and $X_{s}$ the
cost of the repair of the damage per intensity that it causes. Note that
$X_{s}$ is a random variable that depends on where the tornado hits the
ground, say a large city, a medium size town, a village, an open field, etc.
It is assumed that a tornado is equally likely to strike anywhere in the
state. \smallskip

We shall be studying the asymptotic distributional behavior of the randomly
weighted self-normalized L\'{e}vy process $U_{t}/V_{t}$ near 0 and infinity.
Note that $\overline{\Lambda}(0+)=\infty$ implies that $V_{t}>0$ a.s.~for any
$t>0$. Whereas if $\overline{\Lambda}(0+)<\infty$, then, with probability 1,
$V_{t}=0$ for all $t$ close enough to zero. For such $t>0$, $U_{t}/V_{t}$
$=0/0:=0$. Therefore to avoid this triviality, when we consider the asymptotic
behavior of $U_{t}/V_{t}$ near 0 we shall always assume that $\overline
{\Lambda}(0+)=\infty$. \smallskip\ 

Our study is motivated by the following results for weighted sums. Let
$\left\{  Y,Y_{i}:i\geq1\right\}  $ denote a sequence of i.i.d. random
variables, where $Y$ is non-negative and nondegenerate with cdf $G$. Now let
$\left\{  X,X_{i}:i\geq1\right\}  $ be a sequence of i.i.d. random variables,
independent of $\{Y,Y_{i}:i\geq1\}$. Assume that $X$ has cdf $F$ and is in the
class $\mathcal{X}$ of nondegenerate random variables $X$ satisfying
$E|X|<\infty.$ Consider the self-normalized sums
\[
T(n)=\frac{\sum_{i=1}^{n}X_{i}Y_{i}}{\sum_{i=1}^{n}Y_{i}}.
\]
We define $0/0:=0$. Theorem 4 of Breiman \cite{Brei} says that $T(n)$
converges in distribution along the full sequence $\left\{  n\right\}  $ for
\textit{every} $X\in\mathcal{X}$ with at least one limit law being
nondegenerate if and only if $Y\in D(\beta)$, with $0\leq\beta<1$, which means
that for some function $L$ slowly varying at infinity, $P\left\{  Y>y\right\}
=y^{-\beta}L(y),y>0.$ In the case $0<\beta<1$ this is equivalent to $Y\geq0$
being in the domain of attraction of a positive stable law of index $\beta$.
Breiman \cite{Brei} has shown in his Theorem 3 that in this case the limit has
a distribution related to the arcsine law. At the end of his paper Breiman
conjectured that $T(n)$ converges in distribution to a nondegenerate law for
\textit{some} $X\in\mathcal{X}$ if and only if $Y\in D(\beta),$ with
$0\leq\beta<1.$ Mason and Zinn \cite{mz} partially verified his conjecture.
They established the following: \smallskip

Whenever $X$ is nondegenerate and satisfies $E|X|^{p}<\infty$\ for some $p>2,
$ then $T(n)$ converges in distribution to a nondegenerate random variable if
and only if $Y\in D(\beta)$, $0\leq\beta<1$. \smallskip

Recently, Kevei and Mason \cite{KM} investigated the subsequential limits of
$T(n)$. To state their main result we need some definitions. A random variable
$Y$ (not necessarily non-negative) is said to be in the \textit{Feller class}
if there exist sequences of centering and norming constants $\{a_{n}%
\}_{n\geq1}$ and $\{b_{n}\}_{n\geq1}$ such that if $Y_{1},Y_{2},\dots$ are
i.i.d. $Y$ then for every subsequence of $\{n\}$ there exists a further
subsequence $\{n^{\prime}\}$ such that
\[
\frac{1}{b_{n^{\prime}}}\left\{  \sum_{i=1}^{n^{\prime}}Y_{i}-a_{n^{\prime}%
}\right\}  \overset{\mathrm{D}}{\longrightarrow}W,\text{ as }n^{\prime
}\rightarrow\infty,
\]
where $W$ is a nondegenerate random variable. We shall denote this by
$Y\in\mathcal{F}$. Furthermore, $Y$ is in the \textit{centered Feller class},
if $Y$ is in the \textit{Feller class} and one can choose $a_{n}=0$, for all
$n\geq1$. We shall denote this as $Y\in\mathcal{F}_{c}$. The main theorem in
\cite{KM} connects $Y\in\mathcal{F}_{c}$ with the continuity of all of the
subsequential limit laws of $T(n)$. It says that all of the subsequential
distributional limits of $T(n)$ are continuous for any $X$ in the class
$\mathcal{X}$, if and only if $Y\in\mathcal{F}_{c}$.\smallskip

The notions of Feller class and centered Feller class carry over to L\'{e}vy
processes. In particular, a L\'{e}vy process $Y_{t}$ is said to be in the
\textit{Feller class} at infinity if there exists a norming function $B\left(
t\right)  $ and a centering function $A\left(  t\right)  $ such that for each
sequence $t_{k}$ $\rightarrow\infty$ there exists a subsequence $t_{k}%
^{\prime}\rightarrow\infty$ $\ $such that
\[
\left(  Y_{t_{k}^{\prime}}-A(t_{k}^{\prime})\right)  /B(t_{k}^{\prime
})\overset{\mathrm{D}}{\longrightarrow}W,\ \text{ as }k\rightarrow\infty,
\]
where $W$ is a nondegenerate random variable. The L\'{e}vy process $Y_{t}$
belongs to the \textit{centered Feller class} at infinity if\ it is in the
Feller class at infinity and the centering function $A\left(  t\right)  $ can
be chosen to be identically zero. For the definitions of \textit{Feller class}
at zero and \textit{centered Feller class} at zero replace $t_{k}$
$\rightarrow\infty$ and $t_{k}^{\prime}$ $\rightarrow\infty$, by $t_{k}$
$\searrow0$ and $t_{k}^{\prime}$ $\searrow0$, respectively. See Maller and
Mason \cite{MM2} and \cite{MM3} for more details.\medskip

In this paper, we consider the continuous time analog of the results described
above, i.e.~we investigate the asymptotic properties of the self-normalized
L\'{e}vy process
\begin{equation}
T_{t}=U_{t}/V_{t},\label{ratio}%
\end{equation}
as $t\searrow0$ or $t\rightarrow\infty$. The expression \textit{continuous
time analog} is justified\ by Remark 2 in \cite{KM}, where it is pointed out
that under appropriate regularity conditions, norming sequence $\left\{
b_{n}\right\}  _{n\geq1}$ and subsequences $\left\{  n^{\prime}\right\}  $,
\begin{equation}
\left(  \frac{\sum_{1\leq i\leq n^{\prime}t}X_{i}Y_{i}}{b_{n^{\prime}}}%
,\frac{\sum_{1\leq i\leq n^{\prime}t}Y_{i}}{b_{n^{\prime}}}\right)
\overset{\mathrm{D}}{\longrightarrow}(a_{1}t+U_{t},a_{2}t+V_{t}),\text{as
}n^{\prime}\rightarrow\infty.\label{fd}%
\end{equation}
In light of (\ref{fd}) the results that we obtain in the case $t\rightarrow
\infty$ are perhaps not too surprising given those just described for weighted
sums. However, we find our results in the case $t\searrow0$ unexpected.

Our main goal is to establish the following two theorems about the asymptotic
distributional behavior of $U_{t}/V_{t}$. In the process we shall uncover a
lot of information about its subsequential limit laws. First, assuming that
$E|X|^{p}<\infty$, for some $p>2$, we obtain a partial solution to the
continuous time version of the Breiman conjecture, i.e.~the continuous time
version of the result of Mason and Zinn \cite{mz}.

\begin{theorem}
\label{Th2} Assume that $X$ is nondegenerate and for some $p>2$,
$E|X|^{p}<\infty$. Also assume that $\Lambda$ satisfies (\ref{VV}) and, in the
case $t\searrow0$, that $\overline{\Lambda}\left(  0+\right)  =\infty$. The
following are necessary and sufficient conditions for $U_{t}/V_{t}$ to
converge in distribution as $t\searrow0$ (as $t\rightarrow\infty)$ to a random
variable $T$, in which case it must happen that $\left(  EX\right)  ^{2}\leq
ET^{2}\leq EX^{2}$.

(i) $U_{t}/V_{t}\overset{\mathrm{D}}{\rightarrow}T$ and $\left(  EX\right)
^{2}<ET^{2}<EX^{2}$ if and only if $\overline{\Lambda}$ is regularly varying
at zero (infinity) with index $-\beta\in(-1,0)$, in which case the random
variable $T$ has cumulative distribution function
\begin{equation}
P\left\{  T\leq x\right\}  =\frac{1}{2}+\frac{1}{\pi\beta}\arctan\left[
\frac{\int|u-x|^{\beta}\mathrm{sgn}(x-u)F(\mathrm{d}u)}{\int|u-x|^{\beta
}F(\mathrm{d}u)}\tan\frac{\pi\beta}{2}\right]  ,\text{ }x\in\left(
-\infty,\infty\right)  ;\label{integ}%
\end{equation}

(ii) $U_{t}/V_{t}\overset{\mathrm{D}}{\rightarrow}T$ and $ET^{2}=EX^{2}$ if
and only if $\overline{\Lambda}$ is slowly varying at zero (infinity), in
which case $T \overset{\mathrm{D}}{=}X;$

(iii) $U_{t}/V_{t}\overset{\mathrm{D}}{\rightarrow}T$ and $ET^{2}=\left(
EX\right)  ^{2}$ if and only if $\overline{\Lambda}$ is regularly varying at
zero (infinity) with index $-1$, in which case $T=EX.$
\end{theorem}

\begin{remark}
\label{R11}The assumption that $E|X|^{p}<\infty$ for some $p>2$ is only used
in the proof of necessity in Theorem \ref{Th2}. For the sufficiency parts of
the theorem we only need to assume that $X$ is nondegenerate and $E\left\vert
X\right\vert <\infty$. In line with the Breiman \cite{Brei} conjecture we
suspect that $U_{t}/V_{t}\overset{\mathrm{D}}{\rightarrow}T$, as
$t\searrow0$ (as $t\rightarrow\infty$), where $T$ is nondegenerate only if
$\overline{\Lambda}$ satisfies the conclusion of parts (i) or (ii), and in the
case that $T$ is degenerate only if $\overline{\Lambda}$ satisfies the
conclusion of (iii).
\end{remark}

\begin{remark}
\label{R1}A special case of Theorem \ref{Th2} shows that if $W_{t},$ $t>0,$ is
standard Brownian motion, $V_{t}=\inf\left\{  s\geq0:W_{s}>t\right\}  $ and
each $X_{t}$ in (\ref{WL}) is a zero/one random variable $X$ with $P\left\{
X=1\right\}  =1/2$, then$\ U_{t}/V_{t}$ converges in distribution to the
arcsine law as $t\searrow0$ or $t\rightarrow\infty.$ This is a consequence of
the fact that $V_{t}$ is a stable process of index $1/2$, since in this case
we can set $\beta=1/2$ and let $F$ be the cdf of $X$ in (\ref{integ}), which
yields after a little calculation that $T$ has the arcsine density
$g_{T}\left(  t\right)  =\pi^{-1}\left(  t(1-t)\right)  ^{-1/2}$ for $0<t<1$.
Moreover,$\ U_{t}/V_{t}\overset{\mathrm{D}}{=}U_{1}/V_{1}$, for all $t>0$,
which can be seen by using the self-similar property of the $1/2$-stable process.
\end{remark}

\begin{remark}
\label{R1b}Theorem \ref{Th2} has an interesting connection to some results of
Barlow, Pitman and Yor \cite{BPY} and Watanabe \cite{Wat}. Suppose $V_{t}$ is
a strictly stable process of index $0<\beta<1$ and each $X_{t}$ in (\ref{WL})
is a zero/one random variable $X$ with $P\left\{  X=1\right\}  =p$, with
$0<p<1.$ Then $\ $Theorem \ref{Th2} implies that $U_{t}/V_{t}$ converges in
distribution as $t\searrow0$ or $t\rightarrow\infty$ to a random variable
$Y_{\beta,p}$ with density defined for $0<x<1$, by
\[
g_{Y_{\beta,p}}\left(  x\right)  =\frac{\sin\left(  \pi\beta\right)  }{\pi
}\frac{p\left(  1-p\right)  x^{\beta-1}\left(  1-x\right)  ^{\beta-1}}%
{p^{2}\left(  1-x\right)  ^{2\beta}+\left(  1-p\right)  ^{2}x^{2\beta
}+2p\left(  1-p\right)  x^{\beta}\left(  1-x\right)  ^{\beta}\cos\left(
\pi\beta\right)  }.
\]
Furthermore, since $V_{t}$ is self-similar, one sees that $U_{t}/V_{t}%
\overset{\mathrm{D}}{=}U_{1}/V_{1}$, for all $t>0$. Barlow, Pitman and Yor
\cite{BPY} and Watanabe \cite{Wat} show that $g_{Y_{\beta,p}}$ is the density
of the random variable
\[
p^{1/\beta}V_{1}/\left(  p^{1/\beta}V_{1}+\left(  1-p\right)  ^{1/\beta}%
V_{1}^{\prime}\right)  ,
\]
where $V_{1}\overset{\mathrm{D}}{=}V_{1}^{\prime}$ with $V_{1}$ and
$V_{1}^{\prime}$ independent. Moreover, Theorem 2 of Watanabe \cite{Wat} says
that if $A_{t}$ is the occupation time of $Z_{s}$, a $p-$skewed Bessel process
of dimension $2-2\beta$, defined as
\[
A_{t}=\int_{0}^{t}\mathbf{1}\left\{  Z_{s}\geq0\right\}  \mathrm{d}s,
\]
then for all $t>0$, $A_{t}/t$ has a distribution with density $g_{Y_{\beta,p}%
}$. We point out that two additional representations can be given for
$Y_{\beta,p}$ using Propositions \ref{prop1-repr} and \ref{jump-repr} in the
next section. For more about the distribution of $Y_{\beta,p}$ as well as that
of closely related random variables refer to James \cite{Jam}.
\end{remark}

\begin{remark}
\label{R2}Let $V_{t}$ be a subordinator and for each $x\geq0$ let $T\left(
x\right)  $ denote $\inf\left\{  t\geq0:V_{t}>x\right\}  .$ Theorem \ref{Th2}
is analogous to Theorem 6, Chapter 3, of Bertoin \cite{bert}, which says that
$x^{-1}V_{T\left(  x\right)  -}$ converges in distribution as $x\searrow0$,
(as $x\rightarrow\infty$) if and only if $V_{t}$ satisfies the necessary
assumptions of Theorem \ref{Th2} for some $-\beta\in\left[  -1,0\right]  .$
The $\beta=0$ case corresponds to $\overline{\Lambda}$ being slowly varying at
zero (infinity). When $-\beta\in\left(  -1,0\right)  $, the limiting
distribution is the generalized arcsine law.
\end{remark}

Our most significant result about subsequential laws of $U_{t}/V_{t}$ is the
following. Note that contrary to Theorem \ref{Th2} we only assume finite
expectation of $X$.

\begin{theorem}
\label{subseq} Assume $\left(  U_{t},V_{t}\right)  $, $t\geq0$, satisfies
(\ref{E}) and (\ref{VV}) and has joint characteristic function (\ref{holds}).
All subsequential distributional limits of $U_{t}/V_{t}$, as $t\searrow0$, (as
$t\rightarrow\infty$) are continuous for any cdf $F$ in the class $\mathcal{X}
$, if and only if $V_{t}$ is in the centered Feller class at 0 $(\infty)$.
\end{theorem}

\begin{remark}
\label{R4}The proof of Theorem \ref{subseq} shows that if $F$ is in the class
$\mathcal{X}$ and $V_{t}$ is in the centered Feller class at 0 $(\infty)$, all
of the subsequential limit laws of $U_{t}/V_{t}$, as $t\searrow0$, (as
$t\rightarrow\infty$) are not only continuous, but also have Lebesgue
densities on $\mathbb{R}$.
\end{remark}

The rest of the paper is organized as follows. Section 2 contains two
representations of the 2-dimensional L\'{e}vy process $(U_{t},V_{t})$. The
first one plays a crucial role in the proof of Theorem \ref{Th2}, while the
second one points out the connection between the continuous and discrete time
versions of $V_{t}$. We provide a fairly exhaustive list of properties of the
subsequential limit laws of $(U_{t},V_{t})$ in Section 3, and we prove our
main results in Section 4. The Appendix contains some technical results needed
in the proofs.

\section{Preliminaries}

\subsection{Representations for $\left(  U_{t},V_{t}\right)  $}

Let $\left(  U_{t},V_{t}\right)  ,$ $t\geq0$, be a L\'{e}vy process satisfying
(\ref{E}) and (\ref{VV}) with joint characteristic function (\ref{holds}). We
establish two representations for the bivariate L\'{e}vy process.

Let $\varpi_{1},\varpi_{2},\dots$ be a sequence of i.i.d. exponential random
variables with mean $1$, and for each integer $i\geq1$ set $S_{i}=\sum
_{j=1}^{i}\varpi_{j}.$ Independent of $\varpi_{1},\varpi_{2},\dots$ let
$X_{1},X_{2},\dots$ be a sequence of i.i.d. random variables with cdf $F$,
which by (\ref{E}) satisfies $E\left\vert X_{1}\right\vert <\infty.$ Consider
the Poisson process $N(t)$ on $[0,\infty)$ with rate $1$,
\begin{equation}
N\left(  t\right)  =\sum_{j=1}^{\infty}\mathbf{1}_{\{S_{j}\leq t\}}%
\mbox{, }t\geq0.\label{Poisson-proc}%
\end{equation}
Define for $s>0$,
\begin{equation}
\varphi\left(  s\right)  =\sup\left\{  y:\overline{\Lambda}(y)>s\right\}
,\label{phi}%
\end{equation}
where the supremum of the empty set is taken as 0. It is easy to check that
(\ref{VV}) and Lemma \ref{int-trfo} below imply that for all $\delta>0$,
\begin{equation}
\int_{\delta}^{\infty}\varphi\left(  s\right)  \mathrm{d}s<\infty
.\label{delta}%
\end{equation}

We have the following distributional representation of $\left(  U_{t}%
,V_{t}\right)  $:

\begin{proposition}
\label{prop1-repr} For each fixed $t>0$,
\begin{equation}
\left(  U_{t},V_{t}\right)  \overset{\mathrm{D}}{=}\left(  \sum_{i=1}^{\infty
}X_{i}\varphi\left(  \frac{S_{i}}{t}\right)  ,\sum_{i=1}^{\infty}%
\varphi\left(  \frac{S_{i}}{t}\right)  \right)  .\label{ss}%
\end{equation}

\end{proposition}

It is important to note that this representation only holds for fixed $t>0$
and not for the process in $t$. As a first consequence of this representation
we obtain that $E|U_{t}|/V_{t}\leq E|X|<\infty$, in particular, by Markov's
inequality, $U_{t}/V_{t}$ is stochastically bounded.

Now let $\left\{  X_{s}\right\}  _{s\geq0}$ be a class of i.i.d. $F$ random
variables. Consider for each $t\geq0$ the process%
\[
\left(  \sum_{0\leq s\leq t}X_{s}\Delta V_{s},\sum_{0\leq s\leq t}\Delta
V_{s}\right)  ,
\]
where $\Delta V_{s}=V_{s}-V_{s-}$. The following representation reveals the
analogy between the continuous and discrete time self-normalization.

\begin{proposition}
\label{jump-repr} For each fixed $t\geq0$,
\begin{equation}
\left(  U_{t},V_{t}\right)  \overset{\mathrm{D}}{=}\left(  \sum_{0\leq s\leq
t}X_{s}\Delta V_{s},\sum_{0\leq s\leq t}\Delta V_{s}\right)  .\label{RW}%
\end{equation}

\end{proposition}

\begin{remark}
\label{R5}Notice that the process on the right hand side of (\ref{RW}) is a
stationary independent increment process. Since it has the same characteristic
function as $\left(  U_{t},V_{t}\right)  $, the distributional representation
in (\ref{RW}) holds as a process in $t\geq0.$
\end{remark}

\subsection{Proofs of Propositions \ref{prop1-repr} and \ref{jump-repr}}

In the proofs of Propositions \ref{prop1-repr} and \ref{jump-repr} we shall
assume that $\Lambda\left(  \left(  0,\infty\right)  \right)  =\infty$. The
case $\Lambda\left(  \left(  0,\infty\right)  \right)  <\infty$ follows by the
same methods.\medskip

First we state a useful lemma giving a well-known change of variables formula
(see Revuz and Yor \cite{ry}, Prop. 4.9, p.8, or Br\'{e}maud \cite{brem},
p.301), where the integrals are understood to be Riemann--Stieltjes integrals.

\begin{lemma}
\label{int-trfo} Let $h$ be a measurable function defined on $(a,b]$,
$0<a<b<\infty$, and $R$ a measure on $\left(  0,\infty\right)  $ such that
\[
\overline{R}(x):=R\{(x,\infty)\},\ x>0,
\]
is right continuous and $\overline{R}\left(  \infty\right)  =0$. Assume
$\int_{0}^{\infty}\left\vert h\left(  x\right)  \right\vert R\left(
\mathrm{d}x\right)  <\infty$, and define for $s>0$
\[
\varphi\left(  s\right)  = \sup\left\{  y:\overline{R}\left(  y\right)
>s\right\}  ,
\]
where the supremum of the empty set is defined to be $0$. Then we have
\begin{equation}
\int_{0}^{\infty}h\left(  x\right)  R\left(  \mathrm{d}x\right)  =\int
_{0}^{\infty}h\left(  \varphi\left(  s\right)  \right)  \mathrm{d}%
s.\label{chvar}%
\end{equation}

\end{lemma}

\noindent\textbf{Proof of Proposition \ref{prop1-repr}.} We only consider the
process on $[0,1]$.

Applying the L\'{e}vy--It\^{o} integral representation of a L\'{e}vy process
to our case we have that a.s.~for each $t\geq0$
\begin{equation}
(U_{t},V_{t})=\int_{\mathbb{R}^{2}\backslash\{0\}}(u,v)N([0,t],\mathrm{d}%
u,\mathrm{d}v),\label{levy-ito}%
\end{equation}
where $N$ is a Poisson point process on $(0,1)\times\mathbb{R}\times
\lbrack0,\infty)$, with intensity measure $\text{Leb}\times\Pi$, where $\Pi$
is the L\'{e}vy measure as in (\ref{pi}).

For the Poisson point process we have the representation
\begin{equation}
N=\sum_{i=1}^{\infty}\delta_{(U_{i},X_{i}\varphi(S_{i}),\varphi(S_{i}%
))},\label{Pp-repr}%
\end{equation}
where $\{U_{i}\}$ are i.i.d.~Uniform$(0,1)$ random variables, independent from
$\{X_{i}\}$ and $\{\varpi_{i}\}$. (At this step we consider the L\'{e}vy
process on $[0,1].$) To see this, let
\[
M=\sum_{i=1}^{\infty}\delta_{(U_{i},X_{i},S_{i})},
\]
which is a marked Poisson point process on $[0,1]\times\mathbb{R}%
\times(0,\infty)$, with intensity measure $\nu=\text{Leb}\times F\times
\text{Leb}$. Put $h(u,x,s)=(u,x\varphi(s),\varphi(s))$. Then $\nu\circ
h^{-1}=\text{Leb}\times\Pi$. Thus Proposition 2.1 in Rosi\'{n}ski
\cite{Rosinski} implies that the sequences $\{U_{i}\},\{X_{i}\},\{S_{i}\}$ can
be defined on the same space as $N$ such that (\ref{Pp-repr}) holds.

Using (\ref{Pp-repr}) for $N$, from (\ref{levy-ito}) we obtain that a.s.~for
each $t \in[0,1]$
\begin{equation}
\label{proc-repr}(U_{t}, V_{t}) = \sum_{i=1}^{\infty}\left(  X_{i}
\varphi(S_{i}), \varphi(S_{i}) \right)  \mathbf{1}_{ \{ U_{i} \leq t \} }.
\end{equation}

To finish the proof note that if $\sum_{i=1}^{\infty}\delta_{x_{i}}$ is a
Poisson point process and independently $\{\beta_{i}\}$ is an
i.i.d.~Bernoulli$(t)$ sequence, then
\[
\sum_{i=1}^{\infty}\delta_{x_{i}}\mathbf{1}_{\{\beta_{i}=1\}}\overset
{\mathrm{D}}{=}\sum_{i=1}^{\infty}\delta_{x_{i}/t},
\]
i.e.~for a Poisson point process independent Bernoulli thinning and scaling
are distributionally the same.

Since the process representation (\ref{proc-repr}) can be extended to any
finite interval $[0,T]$ (see the final remark in \cite{Rosinski}), this
completes the proof. $\,\hfill\Box$ \medskip

We point out that Proposition \ref{prop1-repr}\textbf{\ }can also be proved by
the same way as Proposition 5.1 in Maller and Mason \cite{MM}.

Next we turn to the proof of the second representation.

\noindent\textbf{Proof of Proposition \ref{jump-repr}.} Let $\left\{
N_{n}\right\}  _{n\geq1}$ be a sequence of independent Poisson processes on
$\left[  0,\infty\right)  $ with rate $1$. Independent of $\left\{
N_{n}\right\}  _{n\geq1}$ let $\left\{  \xi_{i,n}\right\}  _{i\geq1,n\geq1}$
be an array of independent random variables such that for each $i\geq1,n\geq
1$, $\xi_{i,n}$ has distribution $P_{i,n}$ defined for each Borel subset of
$A$ of $\mathbb{R}$ by
\[
P_{i,n}\left(  A\right)  =P\left\{  \xi_{i,n}\in A\right\}  =\Lambda\left(
A\cap\left[  a_{n},a_{n-1}\right)  \right)  /\mu_{n}\text{,}%
\]
where $a_{n}$ is a strictly decreasing sequence of positive numbers converging
to zero such that $a_{0}=\infty$ and for all $n\geq1,$ $0<\mu_{n}%
=\Lambda\left(  \left[  a_{n},a_{n-1}\right)  \right)  <\infty$.

The process $V_{t},$ $t\geq0$, has the representation as the Poisson point
process%
\[
V_{t}=\sum_{n=1}^{\infty}\sum_{i\leq N_{n}\left(  t\mu_{n}\right)  }\xi
_{i,n}=:\sum_{n=1}^{\infty}V_{t}^{\left(  n\right)  }.
\]
See Bertoin \cite{bert}, page 16. In this representation
\[
V_{t}^{\left(  n\right)  }=\sum_{0\leq s\leq t}\Delta V_{s}\mathbf{1}%
_{\{a_{n}\leq\Delta V_{s}<a_{n-1}\}}%
\]
and
\[
\Delta V_{s}\mathbf{1}_{\{a_{n}\leq\Delta V_{s}<a_{n-1}\}}=\sum_{i\leq
N_{n}\left(  s\mu_{n}\right)  }\xi_{i,n}-\sum_{i\leq N_{n}\left(  s\mu
_{n}-\right)  }\xi_{i,n}.
\]
Moreover if $\Delta V_{s}>0$ there exists a unique pair $\left(  i,n\right)  $
such that $\Delta V_{s}=\xi_{i,n}$. Clearly%
\begin{equation}%
\begin{split}
&  \bigg(\sum_{0\leq s\leq t}X_{s}\Delta V_{s}\mathbf{1}_{\{a_{n}\leq\Delta
V_{s}<a_{n-1}\}},\sum_{0\leq s\leq t}\Delta V_{s}\mathbf{1}_{\{a_{n}\leq\Delta
V_{s}<a_{n-1}\}}\bigg)\\
&  \overset{\mathrm{D}}{=}\bigg(\sum_{i\leq N_{n}\left(  t\mu_{n}\right)
}X_{i,n}\,\xi_{i,n},\sum_{i\leq N_{n}\left(  t\mu_{n}\right)  }\xi
_{i,n}\bigg)=:\left(  U_{t}^{\left(  n\right)  },V_{t}^{\left(  n\right)
}\right)  ,
\end{split}
\label{comp}%
\end{equation}
where $\{X_{i,n}\}_{i\geq1,n\geq1}$ is an array of i.i.d. random variables
with common distribution function $F$. Notice that the process $\left(
U_{t}^{\left(  n\right)  },V_{t}^{\left(  n\right)  }\right)  $ in
(\ref{comp}) is a compound Poisson process. Keeping this in mind, we see after
a little calculation that
\[
E\exp\left(  \mathrm{i}\left(  \theta_{1}U_{t}^{\left(  n\right)  }+\theta
_{2}V_{t}^{\left(  n\right)  }\right)  \right)  =\exp\left(  t\int_{\left[
a_{n},a_{n-1}\right)  }\int_{-\infty}^{\infty}\left(  e^{\mathrm{i}(\theta
_{1}u+\theta_{2}v)}-1\right)  F\left(  \mathrm{d}u/v\right)  \Lambda\left(
\mathrm{d}v\right)  \right)  .
\]
Since the random variables $\left\{  \left(  U_{t}^{\left(  n\right)  }%
,V_{t}^{\left(  n\right)  }\right)  \right\}  _{n\geq1}$ are independent we
readily conclude that (\ref{holds}) holds. \hfill$\Box\medskip$

\section{Additional asymptotic distribution results along subsequences}

Let $\mathrm{id}(a,b,\nu)$ denote an infinitely divisible distribution on
$\mathbb{R}^{d}$ with characteristic exponent
\[
\mathrm{i}u^{\prime}b-\frac{1}{2}u^{\prime}au+\int\left(  e^{\mathrm{i}%
u^{\prime}x}-1-\mathrm{i}u^{\prime}x \mathbf{1}_{\{ |x|\leq1 \}} \right)
\nu(\mathrm{d}x),
\]
where $b\in\mathbb{R}^{d}$, $a\in\mathbb{R}^{d\times d}$ is a positive
semidefinite matrix, $\nu$ is a L\'{e}vy measure on $\mathbb{R}^{d}$ and
$u^{\prime}$ stands for the transpose of $u$. In our case $d$ is $1$ or $2$.
For any $h>0$ put
\[
a^{h}=a+\int_{|x|\leq h}xx^{\prime}\nu(\mathrm{d}x)\text{ and }b^{h}%
=b-\int_{h<|x|\leq1}x\nu(\mathrm{d}x).
\]

When $d=1$, id$(b,\Lambda),$ with L\'{e}vy measure $\Lambda$ on $\left(
0,\infty\right)  $, such that (\ref{VV}) holds, and $b\geq0$, denotes a
non-negative infinitely divisible distribution with Laplace transform
\[
\exp\left(  -\theta b-\int_{0}^{\infty}\left(  1-e^{-\theta u}\right)
\Lambda\left(  \mathrm{d}u\right)  \right)  .
\]
In this section it will be convenient to use the following representation for
the joint characteristic function of the L\'{e}vy process $\left(  U_{t}%
,V_{t}\right)  $, $t\geq0,$ satisfying (\ref{E}) and (\ref{VV}) and having
joint characteristic function (\ref{holds}):
\begin{equation}%
\begin{split}
\phi\left(  t,\theta_{1},\theta_{2}\right)   &  =\exp\left(  \mathrm{i}%
t(\theta_{1}b_{1}+\theta_{2}b_{2})\right)  \times\\
&  \exp\left(  t\int_{(0,\infty)}\int_{-\infty}^{\infty}\left(  e^{\mathrm{i}%
(\theta_{1}u+\theta_{2}v)}-1-\left(  \mathrm{i}\theta_{1}u+\mathrm{i}%
\theta_{2}v\right)  \mathbf{1}_{\{u^{2}+v^{2}\leq1\}}\right)  \Pi\left(
\mathrm{d}u,\mathrm{d}v\right)  \right)  ,
\end{split}
\label{cf-uv}%
\end{equation}
where $\Pi\left(  \mathrm{d}u,\mathrm{d}v\right)  $ is as in (\ref{pi}) and
\begin{equation}
\mathbf{b}=\left(
\begin{array}
[c]{c}%
b_{1}\\
b_{2}%
\end{array}
\right)  =\left(
\begin{array}
[c]{c}%
\int_{0<u^{2}+v^{2}\leq1}u\Pi\left(  \mathrm{d}u,\mathrm{d}v\right) \\
\int_{0<u^{2}+v^{2}\leq1}v\Pi\left(  \mathrm{d}u,\mathrm{d}v\right)
\end{array}
\right)  .\label{b}%
\end{equation}
Note that assumptions (\ref{E}) and (\ref{VV}) insure that (\ref{b}) is well
defined.\smallskip

First we investigate the possible subsequential distributional limits of
$\left(  U_{t},V_{t}\right)  $. The following theorem is an analog of Theorem
1 in \cite{KM}.

\begin{theorem}
\label{2-dim-conv} Consider the bivariate L\'{e}vy process $\left(
U_{t},V_{t}\right)  ,$ $t\geq0$, satisfying (\ref{E}) and (\ref{VV}) with
joint characteristic function (\ref{cf-uv}). Assume that for some
deterministic sequences $t_{k}\searrow0$ ($t_{k}\rightarrow\infty)$ and
$B_{k}$ the distributional convergence
\begin{equation}
\frac{V_{t_{k}}}{B_{k}}\overset{\mathrm{D}}{\longrightarrow}V\label{v-conv}%
\end{equation}
holds, where $V$ has $\mathrm{id}(b,\Lambda_{0})$ distribution with L\'{e}vy
measure $\Lambda_{0}$ on $\left(  0,\infty\right)  $. Then
\begin{equation}
\left(  \frac{U_{t_{k}}}{B_{k}},\frac{V_{t_{k}}}{B_{k}}\right)  \overset
{\mathrm{D}}{\longrightarrow}(U,V),\label{uv-conv}%
\end{equation}
where $(U,V)$ has $\mathrm{id}(\mathbf{0},\mathbf{c},\Pi_{0})$ distribution
with L\'{e}vy measure $\Pi_{0}\left(  \mathrm{d}u,\mathrm{d}v\right)
=F\left(  \mathrm{d}u/v\right)  \Lambda_{0}\left(  \mathrm{d}v\right)  $ on
$\mathbb{R}\times\left(  0,\infty\right)  $ and
\begin{equation}
\mathbf{c}=\left(
\begin{array}
[c]{c}%
c_{1}\\
c_{2}%
\end{array}
\right)  =\left(
\begin{array}
[c]{c}%
bEX+\int_{0<u^{2}+v^{2}\leq1}u\Pi_{0}\left(  \mathrm{d}u,\mathrm{d}v\right) \\
b+\int_{0<u^{2}+v^{2}\leq1}v\Pi_{0}\left(  \mathrm{d}u,\mathrm{d}v\right)
\end{array}
\right)  ,\label{c}%
\end{equation}
i.e.~it has characteristic function
\begin{equation}%
\begin{split}
\Psi(\theta_{1},\theta_{2})  &  =Ee^{\mathrm{i}(\theta_{1}U+\theta_{2}V)}%
=\exp\bigg\{\mathrm{i}(\theta_{1}c_{1}+\theta_{2}c_{2})\\
&  \phantom{=}+\int_{0}^{\infty}\int_{-\infty}^{\infty}\left(  e^{\mathrm{i}%
(\theta_{1}u+\theta_{2}v)}-1-\left(  \mathrm{i}\theta_{1}u+\mathrm{i}%
\theta_{2}v\right)  \mathbf{1}_{\{u^{2}+v^{2}\leq1\}}\right)  F\left(
\mathrm{d}u/v\right)  \Lambda_{0}\left(  \mathrm{d}v\right)  \bigg\}.
\end{split}
\label{limit-chfunc-proc}%
\end{equation}

\end{theorem}

Theorem \ref{2-dim-conv} has some immediate consequences concerning the
subsequential limits of $(U_{t},V_{t})$. The first part of the following
corollary is deduced from Theorem \ref{2-dim-conv} and classical theory, i.e.
Theorem 15.14 in \cite{kallenberg}. The second part follows by Fourier inversion.

\begin{corollary}
\label{cor5} Let $(U_{t},V_{t}),$ $t\geq0$, be as in Theorem \ref{2-dim-conv}.
For deterministic constants $t_{k},B_{k}$ the vector $B_{k}^{-1}(U_{t_{k}%
},V_{t_{k}})$ converges in distribution to $\left(  U,V\right)  $ as
$t_{k}\searrow0$ (as $t_{k}\rightarrow\infty)$ having characteristic function
(\ref{limit-chfunc-proc}) if, and only if $t_{k}\overline{\Lambda}%
(vB_{k})\rightarrow\overline{\Lambda}_{0}(v)$ for every continuity point of
$\Lambda_{0}$, and $\int_{0}^{h}xt_{k}\Lambda(\mathrm{d}B_{k}x)\rightarrow
\int_{0}^{h}x\Lambda_{0}(\mathrm{d}x)+b$ for some continuity point $h$ of
$\Lambda_{0}$. Moreover, if $\overline{\Lambda}(0+)=\infty$, or $b>0$ then
$V>0$ a.s., and so $U_{t_{k}}/V_{t_{k}}\overset{}{\overset{\mathrm{D}%
}{\longrightarrow}U/V},$ and with $\Psi$ as in (\ref{limit-chfunc-proc}) for
all $x$
\[
P\left\{  U/V\leq x\right\}  =\frac{1}{2}-\frac{1}{2\pi\mathrm{i}}%
\int_{-\infty}^{\infty}\frac{\Psi(u,-ux)}{u}\mathrm{d}u.
\]

\end{corollary}

The remaining results in this section, though interesting in their own right,
are crucial for the proof of Theorem \ref{subseq}.

The following proposition provides a sufficient condition for $\left(
U,V\right)  $ to have a $C^{\infty}$ 2-dimensional density. It also gives an
alternative proof for Theorem 3 in \cite{KM}. We require the following
notation: Put for $v>0$,%
\begin{equation}
V_{2}\left(  v\right)  =\int_{0<u\leq v}u^{2}\Lambda(\mathrm{d}u).\label{V2}%
\end{equation}

\begin{proposition}
\label{prop-fc0} Assume that $(U,V)$ has joint characteristic function%
\[
Ee^{\mathrm{i}(\theta_{1}U+\theta_{2}V)}=\exp\left\{  \int_{\left(
0,\infty\right)  }\int_{\mathbb{R}}\left(  e^{\mathrm{i}(\theta_{1}%
u+\theta_{2}v)}-1\right)  F\left(  \frac{\mathrm{d}u}{v}\right)
\Lambda(\mathrm{d}v)\right\}  ,
\]
where $\int_{0}^{1}v\Lambda(\mathrm{d}v)<\infty$ and $F$ is in the class
$\mathcal{X}$. Whenever
\begin{equation}
\limsup_{v\searrow0}\frac{v^{2}\overline{\Lambda}\left(  v\right)  }%
{V_{2}\left(  v\right)  }<\infty\label{ss1}%
\end{equation}
holds, then $\left(  U,V\right)  $ has a $C^{\infty}$ density.
\end{proposition}

As a consequence we obtain the following

\begin{corollary}
\label{cor1} Let $(U_{t},V_{t}),$ $t\geq0$, be as in Theorem \ref{2-dim-conv}.
Assume that $V_{t}$ is in the centered Feller class at zero (infinity) and $F$
is in the class $\mathcal{X}$. Then for a suitable norming function $B(t) $
any subsequential distributional limit of
\[
\left(  \frac{U_{t_{k}}}{B(t_{k})},\frac{V_{t_{k}}}{B(t_{k})}\right)
\]
along a subsequence $t_{k}\searrow0$ $(t_{k}\rightarrow\infty)$, say $\left(
W_{1},W_{2}\right)  $, has a $C^{\infty}$ Lebesgue density $f$ on
$\mathbb{R}^{2}$, which implies that the asymptotic distribution of the
corresponding ratio along the subsequence $\{t_{k}\}$ has a Lebesgue density
$g_{T}$ on $\mathbb{R}$.
\end{corollary}

The following corollary is an immediate consequence of Theorem
\ref{2-dim-conv}. Note that a L\'{e}vy process $Y_{t}$ that is in the Feller
class at zero (infinity) but not in the centered Feller class at zero
(infinity) has the required property.

\begin{corollary}
\label{cor2}Let $(U_{t},V_{t}),$ $t\geq0$, be as in Theorem \ref{2-dim-conv}.
Suppose along a subsequence $t_{k}\searrow0$ ($t_{k}\rightarrow\infty)$
\[
\frac{V_{t_{k}}-A(t_{k})}{B(t_{k})}\overset{\mathrm{D}}{\longrightarrow}W,
\]
where $W$ is nondegenerate and $A(t_{k})/B(t_{k})\rightarrow\infty,$ as
$k\rightarrow\infty.$ Then
\[
\frac{U_{t_{k}}}{V_{t_{k}}}\overset{\mathrm{D}}{\longrightarrow}%
EX,\quad\text{as }k\rightarrow\infty.
\]

\end{corollary}

For $t>0$ and $\varepsilon\in(0,1)$ put
\begin{equation}
A_{t}(\varepsilon)=\left\{  \frac{\varphi(S_{1}/t)}{\sum_{i=1}^{\infty}%
\varphi(S_{i}/t)}>1-\varepsilon\right\}  ,\label{AAa}%
\end{equation}
and
\[
\Delta_{t}=\left\vert \frac{\sum_{i=1}^{\infty}X_{i}\varphi(S_{i}/t)}%
{\sum_{i=1}^{\infty}\varphi(S_{i}/t)}-X_{1}\right\vert .
\]

\begin{proposition}
\label{sums-th4} Assume that for a subsequence $t_{k} \searrow0$ or $t_{k}
\to\infty$
\begin{equation}
\label{A-assump}\lim_{\varepsilon\to0} \liminf_{k \to\infty} P \{ A_{t_{k}} (
\varepsilon) \} = \delta> 0,
\end{equation}
then
\[
\lim_{\varepsilon\to0} \liminf_{k \to\infty} P \{ \Delta_{t_{k}}
\leq\varepsilon\} \geq\delta.
\]

\end{proposition}

Together with the stochastic boundedness of $U_{t} / V_{t}$ this implies the following.

\begin{corollary}
\label{cor3}Let $(U_{t},V_{t}),$ $t\geq0$, be as in Theorem \ref{2-dim-conv}.
Assume that (\ref{A-assump}) holds for $V_{t}$, and $P\{X=x_{0}\}>0$ for some
$x_{0}$. Then there exists a subsequence $t_{k}\searrow0$ $(t_{k}%
\rightarrow\infty)$ such that $U_{t_{k}}/V_{t_{k}}\overset{\mathrm{D}%
}{\longrightarrow}T,$ with $P\{T=x_{0}\}>0$.
\end{corollary}

Put
\begin{equation}
R_{t}=\frac{\sum_{i=1}^{\infty}\varphi^{2}\left(  \frac{S_{i}}{t}\right)
}{\left(  \sum_{i=1}^{\infty}\varphi\left(  \frac{S_{i}}{t}\right)  \right)
^{2}}.\label{RT}%
\end{equation}

\begin{proposition}
\label{non-feller} Assume that $R_{t}^{-1}\neq O_{P}(1)$ as $t\searrow0$ or
$t\rightarrow\infty$, then there exists a subsequence $t_{k}\searrow0$ or
$t_{k}\rightarrow\infty$ such that $U_{t_{k}}/V_{t_{k}}\overset{\mathrm{D}%
}{\longrightarrow}T,$ with $P\{T=EX\}>0$.
\end{proposition}

The proofs of Propositions \ref{sums-th4} and \ref{non-feller} are adaptations
of those of Theorems 4 and 5 in \cite{KM}. Therefore we only sketch the proof
of the first one, and omit the proof of the second one.

\section{Proofs of results}

Recall that throughout this paper $\left(  U_{t},V_{t}\right)  $, $t\geq0,$
denotes a L\'{e}vy process satisfying (\ref{E}) and (\ref{VV}) and having
joint characteristic function (\ref{holds}). We start with the proof of
Theorem \ref{2-dim-conv} since this result is crucial for both the proofs of
Theorem \ref{Th2} and \ref{subseq}.

\subsection{Proof of Theorem \ref{2-dim-conv}}

Recall the notation at the beginning of Section 3. Since $V_{t}$ is a
driftless subordinator, by Theorem 15.14 (ii) in \cite{kallenberg},
(\ref{v-conv}) is equivalent to the convergences
\begin{equation}
t_{k}\overline{\Lambda}(vB_{k})\rightarrow\overline{\Lambda}_{0}%
(v),\quad\text{as }k\rightarrow\infty,\label{Levy-conv-V}%
\end{equation}
for any $v>0$ continuity point of $\overline{\Lambda}_{0}$, and
\begin{equation}
\int_{0}^{v}xt_{k}\Lambda(\mathrm{d}B_{k}x)\rightarrow\int_{0}^{v}x\Lambda
_{0}(\mathrm{d}x)+b,\quad\text{as }k\rightarrow\infty,\label{1stmoment-V}%
\end{equation}
where $v>0$ is a fixed continuity point of $\overline{\Lambda}_{0}$.

Notice that using (\ref{cf-uv}) we have that
\[%
\begin{split}
E  &  e^{\mathrm{i}\left(  \theta_{1}\frac{U_{t_{k}}}{B_{k}}+\theta_{2}%
\frac{V_{t_{k}}}{B_{k}}\right)  }=\exp\left\{  \mathrm{i}\frac{t_{k}}{B_{k}%
}(\theta_{1}b_{1}+\theta_{2}b_{2})\right\} \\
&  \phantom{=}\times\exp\left\{  \int\left[  e^{\mathrm{i}(\theta_{1}%
u+\theta_{2}v)/B_{k}}-1-\frac{\mathrm{i}}{B_{k}}(\theta_{1}u+\theta
_{2}v)\mathbf{1}_{\{0<u^{2}+v^{2}\leq1\}}\right]  t_{k}\Pi(\mathrm{d}%
u,\mathrm{d}v)\right\} \\
&  =\exp\left\{  \mathrm{i}\frac{t_{k}}{B_{k}}(\theta_{1}b_{1}+\theta_{2}%
b_{2})\right\} \\
&  \phantom{=}\times\exp\left\{  \int\left[  e^{\mathrm{i}(\theta_{1}%
x+\theta_{2}y)}-1-\mathrm{i}(\theta_{1}x+\theta_{2}y)\mathbf{1}_{\{0<x^{2}%
+y^{2}\leq B_{k}^{-2}\}}\right]  \Pi_{k}(\mathrm{d}x,\mathrm{d}y)\right\}  ,
\end{split}
\]
where $\Pi$ is the L\'{e}vy measure on $\left(  0,\infty\right)
\times\mathbb{R}$ defined by (\ref{pi}) and for each $k\geq1,$ $\Pi_{k}$ is
the L\'{e}vy measure on $\left(  0,\infty\right)  \times\mathbb{R}$ defined by%
\[
\Pi_{k}(\mathrm{d}x,\mathrm{d}y)=t_{k}\Pi(B_{k}\mathrm{d}x,B_{k}\mathrm{d}y).
\]
Further, for each $k\geq0$ and $h>0$ with $\Pi_{0}(\left\{  x:|x|=h\right\}
)=0$, in accordance with the notation at the beginning of Section 3, let
\begin{align*}
a_{k}^{h}  &  =\int_{x^{2}+y^{2}\leq h^{2}}\left(
\begin{array}
[c]{cc}%
x^{2} & xy\\
xy & y^{2}%
\end{array}
\right)  \Pi_{k}(\mathrm{d}x,\mathrm{d}y),\\
b_{k}^{h}  &  =\frac{t_{k}}{B_{k}}\mathbf{b}-\int_{1<x^{2}+y^{2}\leq
B_{k}^{-2}}(x,y)\Pi_{k}(\mathrm{d}x,\mathrm{d}y)-\int_{h^{2}<x^{2}+y^{2}\leq
1}(x,y)\Pi_{k}(\mathrm{d}x,\mathrm{d}y)\text{ }\\
&  =\int_{x^{2}+y^{2}\leq h^{2}}(x,y)\Pi_{k}(\mathrm{d}x,\mathrm{d}y),
\end{align*}
where we used (\ref{b}). We set $a^{h}:=a_{0}^{h}$ and $b^{h}:=b_{0}^{h}.$
\medskip

To show (\ref{uv-conv}), by Theorem 15.14 (i) in \cite{kallenberg} we have to
prove that as $k\rightarrow\infty$,
\begin{equation}
\Pi_{k}\overset{v}{\rightarrow}\Pi_{0}\text{, on }\mathbb{R}^{2}-\left\{
\mathbf{0}\right\} \label{pi-conv-uv}%
\end{equation}
and for some (any) $h>0$ with $\Pi_{0}(\left\{  x:|x|=h\right\}  )=0,$ as
$k\rightarrow\infty$,
\begin{align}
a_{k}^{h}  &  \rightarrow a^{h},\label{a-conv-uv}\\
b_{k}^{h}  &  \rightarrow b^{h}.\label{b-conv-uv}%
\end{align}
To establish (\ref{pi-conv-uv}) it suffices to show that for each $\left(
u,v\right)  $ with $u\geq0$, $v>0$, and $\left(  u,v\right)  $, with $u>0$,
$v=0 $, that when $\left(  u,v\right)  $ is a continuity point of
$\overline{\Pi}_{0}$,
\[
t_{k}\overline{\Pi}(B_{k}u,B_{k}v)\rightarrow\overline{\Pi}_{0}(u,v),\text{ as
}k\rightarrow\infty,\text{ }%
\]
and when $\left(  -u,v\right)  $ is a continuity point of $\Pi_{0}$,
\[
t_{k}\Pi(-B_{k}u,B_{k}v)\rightarrow\Pi_{0}(-u,v),\text{ as }k\rightarrow
\infty;
\]
where for $u\geq0,v>0$,
\[
t_{k}\overline{\Pi}(B_{k}u,B_{k}v)=\int_{v}^{\infty}\overline{F}%
(u/y)t_{k}\Lambda(\mathrm{d}B_{k}y),
\]%
\[
\overline{\Pi}_{0}(u,v)=\int_{v}^{\infty}\overline{F}(u/y)\Lambda
_{0}(\mathrm{d}y),
\]%
\[
t_{k}\Pi(-B_{k}u,B_{k}v)=\int_{v}^{\infty}F(-u/y)t_{k}\Lambda(\mathrm{d}%
B_{k}y)
\]
and
\[
\Pi_{0}(-u,v)=\int_{v}^{\infty}F(-u/y)\Lambda_{0}(\mathrm{d}y).
\]
This follows with obvious changes of notation exactly as in the proof of
Proposition 1 in \cite{KM}.

The proofs that (\ref{a-conv-uv}) and (\ref{b-conv-uv}) hold follow exactly as
in Propositions 2 and 3 in \cite{KM}. It turns out that $a^{h}$ converges to
the zero matrix as $h\searrow0$ and by (\ref{1stmoment-V})
\[
b^{h}=\left(
\begin{array}
[c]{c}%
bEX+\int_{0}^{h}\psi(v)\Lambda_{0}(\mathrm{d}v)\\
b+\int_{0}^{h}\phi(v)v\Lambda_{0}(\mathrm{d}v)
\end{array}
\right)  ,
\]
where $\psi$ and $\phi$ are the following functions of $v\in(0,h]$:
\[
\phi\left(  v\right)  =\int_{\left[  -\sqrt{h^{2}-v^{2}},\sqrt{h^{2}-v^{2}%
}\right]  }F\left(  \frac{\mathrm{d}u}{v}\right)  \text{ and }\psi\left(
v\right)  =\int_{\left[  -\sqrt{h^{2}-v^{2}},\sqrt{h^{2}-v^{2}}\right]
}uF\left(  \frac{\mathrm{d}u}{v}\right)  .
\]
(Refer to \cite{KM} for details.) Thus
\[
\lim_{h\rightarrow0}b^{h}=\left(
\begin{array}
[c]{c}%
bEX\\
b
\end{array}
\right)  ,
\]
and the theorem follows with the stated constants. $\,$ \hfill$\Box$ \medskip

\subsection{Proof of Theorem \ref{Th2}}

The following three lemmas establish the \textquotedblleft in
which\ case\textquotedblright\ parts of (i), (ii) and (iii) of Theorem
\ref{Th2}.

\begin{lemma}
\label{cor4} If $\overline{\Lambda}$ is regularly varying at zero (infinity)
with index $-\beta$ with $\beta\in\left(  0,1\right)  $, then for an
appropriate norming function $B_{t}$ the random variable $B_{t}^{-1}%
(U_{t},V_{t})$ converges in distribution as $t\searrow0$ (as $t\rightarrow
\infty$) to $\left(  U,V\right)  $, having joint characteristic function
\begin{equation}
\phi\left(  \theta_{1},\theta_{2}\right)  =\exp\left(  \int_{(0,\infty)}%
\int_{-\infty}^{\infty}\left(  e^{\mathrm{i}(\theta_{1}u+\theta_{2}%
v)}-1\right)  F\left(  \mathrm{d}u/v\right)  \beta v^{-1-\beta}\mathrm{d}%
v\right) \label{stable}%
\end{equation}
and thus
\begin{equation}
T_{t}=\frac{U_{t}}{V_{t}}\overset{\mathrm{D}}{\longrightarrow}\frac{U}%
{V}\text{, as }t\searrow0\ (\text{as }t\rightarrow\infty).\label{UV}%
\end{equation}
Moreover, the cdf of $U/V$ is given by (\ref{integ}).
\end{lemma}

\noindent\textbf{Proof.} We can find a function $B_{t}$ on $\left[
0,\infty\right)  $ such that
\[
B_{t}=L^{\ast}\left(  t\right)  t^{1/\beta},\text{ }t>0\text{,}%
\]
with $L^{\ast}$ defined on $\left[  0,\infty\right)  $ slowly varying at zero
(infinity) satisfying for all $y>0$,
\[
\overline{\mu}_{t}\left(  y\right)  :=t\overline{\Lambda}\left(
yB_{t}\right)  \rightarrow\overline{\Lambda}_{0}\left(  y\right)  =y^{-\beta
}\text{, as }t\searrow0\text{ (as }t\rightarrow\infty\text{).}%
\]
It is routine to show using well-known properties of regularly varying
functions that for any $y>0$, as $t\searrow0$ (as $t\rightarrow\infty$)
\[
a_{t}^{h}:=\int_{0<y\leq h}y\mu_{t}\left(  \mathrm{d}y\right)  \rightarrow
\frac{\beta h^{1-\beta}}{1-\beta}=\int_{0<y\leq h}y\Lambda_{0}\left(
\mathrm{d}y\right)  =:a^{h}.
\]
Thus by applying Theorem 15.14 (ii) in \cite{kallenberg} we find that
$B_{t}^{-1}V_{t}$ converges in distribution as $t\searrow0$ (as $t\rightarrow
\infty$) to $V$, having characteristic function $\phi\left(  0,\theta
_{2}\right)  .$ This says that $V$ is a subordinator with an id$(0,\Lambda
_{0})$ distribution. Theorem \ref{2-dim-conv} completes the proof of
(\ref{stable}).\smallskip

Next, using Fubini's theorem and the explicit formula for the $\beta$-stable
characteristic function (Meerschaert and Scheffler \cite{meerschaert} p.266),
we have for an appropriate constant $c>0$
\begin{align*}
&  \int_{(0,\infty)}\int_{-\infty}^{\infty}\left(  e^{\mathrm{i}(\theta
_{1}u+\theta_{2}v)}-1\right)  F\left(  \mathrm{d}u/v\right)  \beta
v^{-1-\beta}\mathrm{d}v\\
&  =\int_{-\infty}^{\infty} F(\mathrm{d} u)\int_{0}^{\infty}\left[
e^{\mathrm{i}(\theta_{1}u+\theta_{2})y}-1\right]  \Lambda_{0}(\mathrm{d}y)\\
&  =-c\int_{-\infty}^{\infty}|\theta_{1}u+\theta_{2}|^{\beta}\left(
1-\mathrm{i}\,\mathrm{sgn}(\theta_{1}u+\theta_{2})\tan\frac{\pi\beta}%
{2}\right)  F(\mathrm{d}u).
\end{align*}

We see now that the characteristic function of $U-Vx$ is
\begin{equation}
\label{f(t,x)}Ee^{\mathrm{i}t(U-Vx)} = \exp\left\{  -|t|^{\beta}%
c\int|u-x|^{\beta}F(\mathrm{d}u)\left[  1-\mathrm{i}\,\mathrm{sgn}\left(
t\right)  \tan\frac{\pi\beta}{2}\frac{\int|u-x|^{\beta}\mathrm{sgn}%
(u-x)F(\mathrm{d}u)}{\int|u-x|^{\beta}F(\mathrm{d}u)}\right]  \right\}  .
\end{equation}

Proposition 4 in \cite{Brei} now shows that $T$ has cdf (\ref{integ}).

\hfill$\Box\medskip$

\begin{lemma}
\label{sv} If $\overline{\Lambda}$ is slowly varying at zero (at infinity),
then
\begin{equation}
T_{t}=\frac{U_{t}}{V_{t}}\overset{\mathrm{D}}{\longrightarrow}X\text{, as
}t\searrow0\ (\text{as }t\rightarrow\infty),\label{AA}%
\end{equation}
where in the $t\searrow0$ case we also assume $\overline{\Lambda}(0+)=\infty$.
\end{lemma}

\noindent\textbf{Proof.} The proof follows the lines of that of Lemma 5.3 in
\cite{MM}.

We shall only prove the $t\rightarrow\infty$ case. The $t\searrow0$ case is
nearly identical. Now $\overline{\Lambda}$ slowly varying at infinity implies
that $\varphi$ is non-increasing and rapidly varying at $0$ with index
$-\infty$. (See the argument in Item 5 on p.22 of de Haan \cite{haan}.) This
means that for all $0<\lambda<1$
\[
\varphi\left(  x\lambda\right)  /\varphi\left(  x\right)  \rightarrow
\infty,\mbox{ as }x\searrow0.
\]
By Theorem 1.2.1 on p. 15 of \cite{haan}, rapidly varying at $0$ with index
$-\infty$ implies that
\begin{equation}
\frac{\int_{x}^{\overline{\Lambda}(0+)}\varphi\left(  y\right)  \mathrm{d}%
y}{x\varphi\left(  x\right)  }\rightarrow0,\mbox{ as }x\searrow0.\label{x}%
\end{equation}
By Lemma \ref{psi-lemma} in the Appendix, we have
\begin{align*}
E\left(  \frac{\sum_{i=2}^{\infty}\left\vert X_{i}\right\vert \varphi\left(
\frac{S_{i}}{t}\right)  }{\varphi\left(  \frac{S_{1}}{t}\right)  }%
\Bigg|S_{1}\right)   &  =E\left\vert X\right\vert E\left(  \frac{\sum
_{i=2}^{\infty}\varphi\left(  \frac{S_{i}}{t}\right)  }{\varphi\left(
\frac{S_{1}}{t}\right)  }\Bigg|S_{1}\right) \\
&  =E\left\vert X\right\vert S_{1}\frac{\int_{S_{1}/t}^{\overline{\Lambda
}(0+)}\varphi\left(  y\right)  \mathrm{d}y}{\frac{S_{1}}{t}\varphi\left(
\frac{S_{1}}{t}\right)  },
\end{align*}
and by (\ref{x})
\[
E\left\vert X\right\vert S_{1}\frac{\int_{S_{1}/t}^{\overline{\Lambda}%
(0+)}\varphi\left(  y\right)  \mathrm{d}y}{\frac{S_{1}}{t}\varphi\left(
\frac{S_{1}}{t}\right)  }\overset{\mathrm{P}}{\rightarrow}0,\mbox{ as
}t\rightarrow\infty.
\]
From this we can readily conclude that
\begin{equation}
\sum_{i=1}^{\infty}\varphi\left(  \frac{S_{i}}{t}\right)  =\varphi\left(
\frac{S_{1}}{t}\right)  \left(  1+o_{P}\left(  1\right)  \right)  ,\mbox{ as
}t\rightarrow\infty,\label{lo1}%
\end{equation}
and
\begin{equation}
\sum_{i=1}^{\infty}X_{i}\varphi\left(  \frac{S_{i}}{t}\right)  =X_{1}%
\varphi\left(  \frac{S_{1}}{t}\right)  \left(  1+o_{P}\left(  1\right)
\right)  ,\mbox{ as
}t\rightarrow\infty.\label{lo2}%
\end{equation}
From the representation (\ref{ss}), (\ref{lo1}) and (\ref{lo2}) we see that%
\[
\frac{U_{t}}{V_{t}}\overset{D}{=}\frac{X_{1}\varphi\left(  \frac{S_{1}}%
{t}\right)  \left(  1+o_{P}\left(  1\right)  \right)  }{\varphi\left(
\frac{S_{1}}{t}\right)  \left(  1+o_{P}\left(  1\right)  \right)  }%
=X_{1}+o_{P}\left(  1\right)  ,\mbox{ as }t\rightarrow\infty.
\]
Obviously $T_{t}$ converges in distribution as $t\rightarrow\infty$ to $X$.
\hfill$\Box$

\begin{lemma}
\label{beta11} If $\overline{\Lambda}$ is regularly varying at zero (at
infinity) with index $-1$,
\begin{equation}
T_{t}=\frac{U_{t}}{V_{t}}\overset{\mathrm{D}}{\longrightarrow}EX\text{, as
}t\searrow0\ (\text{as }t\rightarrow\infty).\label{EX}%
\end{equation}

\end{lemma}

\noindent\textbf{Proof.} Since $\overline{\Lambda}$ is regularly varying at
zero (at infinity) with index $-1$, we can find norming and centering
functions $b\left(  t\right)  $ and $a\left(  t\right)  $ such that
$b(t)/a(t)\rightarrow0$ as $t\searrow0\ ($as $t\rightarrow\infty)$ and
\[
b\left(  t\right)  ^{-1}\left(  V_{t}-a\left(  t\right)  \right)
\overset{\mathrm{D}}{\longrightarrow}W,\text{ as }t\searrow0\ (\text{as
}t\rightarrow\infty).
\]
(Here we apply part (i) of Theorem 15.14 in \cite{kallenberg}.) From this we
see that
\[
V_{t}/a(t)\overset{\mathrm{P}}{\longrightarrow}1,\text{ as }t\searrow
0\ (\text{as }t\rightarrow\infty).
\]
A straightforward application of Theorem \ref{2-dim-conv} now shows that%
\[
\left(  \frac{U_{t}}{a(t)},\frac{V_{t}}{a(t)}\right)  \overset{\mathrm{P}%
}{\longrightarrow}\left(  EX,1\right)  \text{, as }t\searrow0\ (\text{as
}t\rightarrow\infty).
\]
\hfill$\Box$

Next we turn to the necessary and sufficient parts of (i), (ii) and (iii).
Assume that for some random variable $T$
\begin{equation}
T_{t}\overset{\mathrm{D}}{\longrightarrow}T,\text{ as }t\searrow0\ (\text{as
}t\rightarrow\infty),\label{conT}%
\end{equation}
where in the case $t\searrow0$ we assume that $\overline{\Lambda}(0+)=\infty$.
Our basic tool will be Proposition \ref{prop1-repr}, which says that
\begin{equation}
T_{t}=\frac{U_{t}}{V_{t}}\overset{\mathrm{D}}{=}\frac{\sum_{i=1}^{\infty}%
X_{i}\varphi\left(  \frac{S_{i}}{t}\right)  }{\sum_{i=1}^{\infty}%
\varphi\left(  \frac{S_{i}}{t}\right)  }.\label{repr}%
\end{equation}
Since we assume that
\begin{equation}
E\left\vert X\right\vert ^{p}<\infty\label{delta1}%
\end{equation}
for some $p>2$, we get by Jensen's inequality that
\[
E\left\vert T_{t}\right\vert ^{p}\leq E\left\vert X\right\vert ^{p}<\infty.
\]
(This is the only place in the proof that we use assumption (\ref{delta1}).)
Notice that (\ref{conT}) and (\ref{delta1}) imply that
\begin{equation}
ET_{t}^{2}\rightarrow ET^{2},\text{ as }t\searrow0\ (\text{as }t\rightarrow
\infty).\label{varconv}%
\end{equation}
Obviously $ET_{t}=EX$ and a little calculation gives that
\[
ET_{t}^{2}=(EX)^{2}+{Var}(X)ER_{t},
\]
where $R_{t}$ is defined as in (\ref{RT}). Clearly, $R_{t}\in\lbrack0,1]$ and
whenever (\ref{varconv}) holds, then for some $0\leq\beta\leq1$
\begin{equation}
ER_{t}\rightarrow1-\beta\text{, as }t\searrow0\text{ (as }t\rightarrow
\infty),\label{rt}%
\end{equation}
which is equivalent to
\begin{equation}
\left(  EX\right)  ^{2}\leq ET^{2}\leq EX^{2}.\label{vE}%
\end{equation}

It turns out that the value of $0\leq\beta\leq1$ determines the asymptotic
distribution of $T_{t}$ as $t\searrow0$ (as $t\rightarrow\infty)$ and the
behavior of the L\'{e}vy function $\overline{\Lambda}$ near zero (at
infinity). For instance, when $\beta=1$, $Var\left(  T_{t}\right)
\rightarrow0 $, which implies that
\begin{equation}
T_{t}\overset{\mathrm{P}}{\longrightarrow}EX\text{, as }t\searrow0\text{ (as
}t\rightarrow\infty).\label{beta1}%
\end{equation}
In general we have the following result, which in combination with Lemmas
\ref{cor4}, \ref{sv} and \ref{beta11} will complete the proof of Theorem
\ref{Th2}.

\begin{proposition}
\label{prop-rt} If (\ref{rt}) holds for some $0\leq\beta\leq1$, then
$\overline{\Lambda}$ is regularly varying at zero (infinity) with index
$-\beta$. (In the case $t\searrow0$ we assume $\overline{\Lambda}(0+)=\infty$.)
\end{proposition}

\noindent\textbf{Proof.} Recall the definition of $N(t)$ in
(\ref{Poisson-proc}) and notice that by (\ref{RT}) for any $t>0$ we can write
\[
R_{t}=\frac{\int_{0}^{\infty}\varphi^{2}\left(  s\right)  N(\mathrm{d}
ts)}{\left(  \int_{0}^{\infty}\varphi\left(  s\right)  N(\mathrm{d}
ts)\right)  ^{2}}.
\]
Define for $T>0$ its truncated version
\begin{equation}
R_{t}(T)=\frac{\int_{0}^{T}\varphi^{2}\left(  s\right)  N( \mathrm{d}
ts)}{\left(  \int_{0}^{T}\varphi\left(  s\right)  N( \mathrm{d} ts)\right)
^{2}}.\label{RT1}%
\end{equation}
Given that $N(Tt)=n$
\[
R_{t}(T)\overset{\mathrm{D}}{=}\frac{\sum_{i=1}^{n}\varphi^{2}(V_{i})}{\left(
\sum_{i=1}^{n}\varphi(V_{i})\right)  ^{2}},
\]
where $V_{1},\ldots,V_{n}$ are i.i.d. Uniform$(0,T)$. The same computation as
in Maller and Mason \cite{MM} gives%
\[
ER_{t}(T)=t\int_{0}^{\infty}u\left(  \int_{0}^{T}\varphi^{2}(s)e^{-u\varphi
(s)}\mathrm{d}s\right)  e^{-t\int_{0}^{T}(1-e^{-u\varphi(s)})\mathrm{d}%
s}\,\mathrm{d}u.
\]
Clearly $R_{t}(T)\leq1$. Also $R_{t}(T)\overset{\mathrm{D}}{\rightarrow}R_{t}
$ as $T\rightarrow\infty$ and thus
\begin{equation}
ER_{t}(T)\rightarrow ER_{t},\text{ as }T\rightarrow\infty.\label{conv}%
\end{equation}
For each $T>0$ and $u>0$, set%
\[
\Phi_{T}\left(  u\right)  =\int_{0}^{T}(1-e^{-u\varphi(s)})\mathrm{d}s\text{,
}\Phi\left(  u\right)  =\int_{0}^{\infty}(1-e^{-u\varphi(s)})\mathrm{d}s\text{
and}%
\]%
\begin{equation}
f_{T,t}\left(  u\right)  =-tu\Phi_{T}^{\prime\prime}\left(  u\right)
e^{-t\Phi_{T}\left(  u\right)  }=tu\left(  \int_{0}^{T}\varphi^{2}%
(s)e^{-u\varphi(s)}\mathrm{d}s\right)  e^{-t\int_{0}^{T}(1-e^{-u\varphi
(s)})\mathrm{d}s}.\label{fT}%
\end{equation}
Also for $u>0$, set
\begin{equation}
f_{\left(  t\right)  }\left(  u\right)  =-tu\Phi^{\prime\prime}\left(
u\right)  e^{-t\Phi\left(  u\right)  }=tu\left(  \int_{0}^{\infty}\varphi
^{2}(s)e^{-u\varphi(s)}\mathrm{d}s\right)  e^{-t\int_{0}^{\infty
}(1-e^{-u\varphi(s)})\mathrm{d}s}.\label{def-f(u)}%
\end{equation}
We have in this notation,
\begin{equation}
ER_{t}(T)=\int_{0}^{\infty}f_{T,t}\left(  u\right)  \mathrm{d}u.\label{ERT}%
\end{equation}

\noindent\textit{Case 1: $\beta\in\lbrack0,1)$}. In this case we must first
show that as $T\rightarrow\infty$
\begin{equation}
ER_{t}(T)=\int_{0}^{\infty}f_{T,t}\left(  u\right)  \mathrm{d}u\rightarrow
\int_{0}^{\infty}f_{\left(  t\right)  }\left(  u\right)  \mathrm{d}%
u,\label{PL}%
\end{equation}
which by (\ref{conv}) implies
\begin{equation}
\int_{0}^{\infty}f_{\left(  t\right)  }\left(  u\right)  \mathrm{d}%
u=ER_{t}.\label{int1}%
\end{equation}
It turns out to be surprisingly tricky to justify the passing-to-the-limit in
(\ref{PL}). Lemma \ref{lemma-rt-inf} and Proposition \ref{P1} in the Appendix
handle this problem, and imply that expression (\ref{int1}) is valid for
$ER_{t}$. After this identity is established, the proof is completed by an
easy modification of that of Proposition 5.2 in \cite{MM}, which is based on
Tauberian theorems. Therefore we omit it. \smallskip

\noindent\textit{Case 2: $\beta=1$}. In this case, it is not necessary to
verify (\ref{PL}). Note that we have that by (\ref{rt}) with $\beta=1$
\[
ER_{t}\rightarrow0\text{, as }t\searrow0\text{ }(t\rightarrow\infty).
\]
Therefore since
\[
ER_{t}(T)\rightarrow ER_{t}\geq\int_{0}^{\infty}f_{\left(  t\right)  }\left(
u\right)  \mathrm{d}u,
\]
we can conclude that as $t\searrow0$ ($t\rightarrow\infty$),
\begin{equation}
-t\int_{0}^{\infty}u\Phi^{\prime\prime}(u)e^{-t\Phi(u)}\mathrm{d}u=\int
_{0}^{\infty}f_{\left(  t\right)  }\left(  u\right)  \mathrm{d}u\rightarrow
0,\label{z}%
\end{equation}
which is all we need for the following argument to work for $\beta=1.$
Applying Lemma \ref{int-trfo}, we get
\[
\Phi(u)=\int_{0}^{\infty}(1-e^{-ux})\Lambda\left(  \mathrm{d}x\right)  ,
\]
which by integrating by parts and using (\ref{VV}) is equal to
\[
\Phi(u)=u\int_{0}^{\infty}\overline{\Lambda}(y)e^{-uy}\mathrm{d}y.
\]
Let $q(y)$ denote the inverse function of $\Phi$. From the expression for
$f_{\left(  t\right)  }\left(  u\right)  $ in (\ref{def-f(u)}) and (\ref{z})
we obtain
\[
t^{-1}\int_{0}^{\infty}f_{\left(  t\right)  }\left(  u\right)  \mathrm{d}%
u=-\int_{0}^{\infty}e^{-ty}q(y)\Phi^{\prime\prime}(q(y))q(\mathrm{d}y)\sim
o\left(  t^{-1}\right)  ,
\]
as $t\rightarrow0$ ($t\rightarrow\infty$). Using Theorem 1.7.1 (Theorem
1.7.1') in Bingham et al \cite{bgt} we obtain
\[
-\int_{0}^{x}q(y)\Phi^{\prime\prime}(q(y))q(\mathrm{d}y)\sim o\left(
x\right)  ,
\]
as $x\rightarrow\infty$ ($x\rightarrow0$). Changing the variables and putting
$x=\Phi(v)$ we have
\[
-\int_{0}^{v}u\Phi^{\prime\prime}(u)\mathrm{d}u=o\left(  \Phi(v)\right)  ,
\]
as $v\rightarrow\infty$ ($v\rightarrow0$). Integrating by parts we get
\[
-\int_{0}^{v}u\Phi^{\prime\prime}(u)\mathrm{d}u=-v\Phi^{\prime}(v)+\Phi
(v)=o\left(  \Phi(v)\right)  ,
\]
which gives
\[
\frac{v\Phi^{\prime}(v)}{\Phi(v)}\rightarrow1,
\]
as $v\rightarrow\infty$ ($v\rightarrow0$). This last limit readily implies
that
\[
v^{-1}\Phi(v)=\int_{0}^{\infty}\overline{\Lambda}(y)e^{-vy}\mathrm{d}y
\]
is slowly varying at infinity (zero). By Theorem 1.7.1' (Theorem 1.7.1) in
\cite{bgt} we obtain that $\int_{0}^{x}\overline{\Lambda}(y)\mathrm{d}y$ is
slowly varying at zero (infinity), which by Theorem 1.7.2.b (Theorem 1.7.2) in
\cite{bgt} implies that $\overline{\Lambda}$ is regularly varying at zero with
index $-1$ (at infinity). \hfill$\Box$

\subsection{Proof of Theorem \ref{subseq}}

Before we proceed with the proofs it will be helpful to first cite some
results from Maller and Mason \cite{MM2}, \cite{MM3} and \cite{MM4}.\medskip

Let $Y_{t}$ be a L\'{e}vy process with L\'{e}vy triplet $\ (\sigma^{2}%
,\gamma,\nu)$, i.e.~$Y_{1}$ has $\mathrm{id}(\sigma^{2},\gamma,\nu)$
distribution$.$ Theorem 1 in Maller and Mason \cite{MM2} states $Y_{t}$
belongs to the Feller class at infinity, if and only if
\begin{equation}
\limsup_{x\rightarrow\infty}\frac{x^{2}\nu\{(-\infty,-x)\cup(x,\infty
)\}}{\sigma^{2}+\int_{|y|\leq x}y^{2}\nu(\mathrm{d}y)}<\infty,\label{fcc}%
\end{equation}
and furthermore $Y_{t}$ belongs to the \textit{centered Feller class} at
infinity if and only if
\begin{equation}
\limsup_{x\rightarrow\infty}\frac{x^{2}\nu\{(-\infty,-x)\cup(x,\infty
)\}+x\left\vert \gamma+\int_{1<|y|\leq x}y\nu(\mathrm{d}y)\right\vert }%
{\sigma^{2}+\int_{|y|\leq x}y^{2}\nu(\mathrm{d}y)}<\infty.\label{cfc}%
\end{equation}
For the corresponding equivalences of \textit{Feller class} at zero and
\textit{centered Feller class} at zero replace $x\rightarrow\infty$ by
$x\searrow0$, respectively; see Theorems 2.1 and 2.3 in \cite{MM3}.

It turns out by using the assumption that $V_{t}$ is a subordinator and by
arguing as in the proof of Propositions \ref{prop1-repr} or of Proposition 5.1
in \cite{MM} we get that
\[
\sqrt{R_{t}^{-1}}=\frac{\sum_{i=1}^{\infty}\varphi\left(  \frac{S_{i}}%
{t}\right)  }{\sqrt{\sum_{i=1}^{\infty}\varphi^{2}\left(  \frac{S_{i}}%
{t}\right)  }}\overset{\mathrm{D}}{=}\frac{V_{t}}{\sqrt{\sum_{0\leq s\leq
t}\left(  \Delta V_{t}\right)  ^{2}}}.
\]
From this distributional equality one can conclude that $\sqrt{R_{t}^{-1}}$ is
stochastically bounded as $t\searrow0$ $\left(  t\rightarrow\infty\right)  $
if and only if
\begin{equation}
\limsup_{t\searrow0\text{ }\left(  t\rightarrow\infty\right)  }\frac{t\int
_{0}^{t}x\Lambda(\mathrm{d}x)}{\int_{0}^{t}x^{2}\Lambda(\mathrm{d}%
x)+t^{2}\overline{\Lambda}(t)}<\infty.\label{sc-condition}%
\end{equation}
by applying Theorem 3.1 in \cite{MM4} in the case $t\rightarrow\infty$, and
Proposition 5.1 in \cite{MM3} (with $a(t)\equiv0$ there, and a small
modification) when $t\searrow0$. The partial sum version of this result was
proved by Griffin \cite{griff}.\medskip

\textbf{Proof of Proposition \ref{prop-fc0}.} We first assume $X$ is
nondegenerate and $EX=0$, which implies that there is an $a\geq1$ such that
\begin{equation}
F\left(  a\right)  -F\left(  0\right)  >0\text{ and }F\left(  0\right)
-F\left(  -a\right)  >0.\label{aa}%
\end{equation}
We need the following lemma.

\begin{lemma}
Whenever (\ref{ss1}) holds and $X$ is nondegenerate and $EX=0$, there exist
$0<\kappa<1$ and $d>0$ such that with $a\geq1$ as in (\ref{aa}), if $2a\left(
\left\vert \theta_{1}\right\vert \vee\left\vert \theta_{2}\right\vert \right)
\geq1$, then
\begin{equation}
\mathfrak{Re}\left\{  \int_{\left(  0,\infty\right)  }\int_{\mathbb{R}}\left(
e^{\mathrm{i}(\theta_{1}x+\theta_{2}v)}-1\right)  F\left(  \frac{\mathrm{d}%
x}{v}\right)  \Lambda(\mathrm{d}v)\right\}  \leq-d\left(  \left\vert
\theta_{1}\right\vert ^{\kappa}+\left\vert \theta_{2}\right\vert ^{\kappa
}\right)  .\label{ss3}%
\end{equation}

\end{lemma}

\noindent\textbf{Proof.} Notice that
\[
\mathfrak{Re}\int_{\left(  0,\infty\right)  }\int_{\mathbb{R}}\left(
e^{\mathrm{i}(\theta_{1}x+\theta_{2}v)}-1\right)  F\left(  \frac{\mathrm{d}%
x}{v}\right)  \Lambda(\mathrm{d}v)=\int_{\left(  0,\infty\right)  }%
\int_{\mathbb{R}}\left(  \cos(\theta_{1}x+\theta_{2}v)-1\right)  F\left(
\frac{\mathrm{d}x}{v}\right)  \Lambda(\mathrm{d}v)
\]%
\[
\leq\int_{0\leq v\leq1/\left(  2a\left(  \left\vert \theta_{1}\right\vert
\vee\left\vert \theta_{2}\right\vert \right)  \right)  }\int_{\left\vert
x\right\vert \leq va}\left(  \cos(\theta_{1}x+\theta_{2}v)-1\right)  F\left(
\frac{\mathrm{d}x}{v}\right)  \Lambda(\mathrm{d}v).
\]
Observe that whenever $\left\vert x\right\vert \leq av$ with $a\geq1$ and
$0\leq v\leq1/\left(  2a\left(  \left\vert \theta_{1}\right\vert
\vee\left\vert \theta_{2}\right\vert \right)  \right)  $,
\[
\left\vert \theta_{1}x\right\vert +\left\vert \theta_{2}v\right\vert
\leq\left(  \left\vert a\theta_{1}\right\vert +\left\vert \theta
_{2}\right\vert \right)  v\leq1.
\]
For some $c>0$,
\[
\sup_{0\leq\left\vert u\right\vert \leq1}\frac{\cos u-1}{u^{2}}\leq-c,
\]
thus
\[
\int_{0\leq v\leq1/\left(  2a\left(  \left\vert \theta_{1}\right\vert
\vee\left\vert \theta_{2}\right\vert \right)  \right)  }\int_{\left\vert
x\right\vert \leq va}\left(  \cos(\theta_{1}x+\theta_{2}v)-1\right)  F\left(
\frac{\mathrm{d}x}{v}\right)  \Lambda(\mathrm{d}v)
\]%
\[
\leq-c\int_{0\leq v\leq1/\left(  2a\left(  \left\vert \theta_{1}\right\vert
\vee\left\vert \theta_{2}\right\vert \right)  \right)  }\int_{\left\vert
x\right\vert \leq av}\left(  \theta_{1}x+\theta_{2}v\right)  ^{2}F\left(
\frac{\mathrm{d}x}{v}\right)  \Lambda(\mathrm{d}v).
\]
Now when $\theta_{1}\theta_{2}\geq0$ we have $\theta_{1}\theta_{2}\int_{0\leq
x\leq va}xF\left(  \frac{\mathrm{d}x}{v}\right)  \geq0$, and we get that the
last bound is
\[
\leq-c\int_{0\leq v\leq1/\left(  2a\left(  \left\vert \theta_{1}\right\vert
\vee\left\vert \theta_{2}\right\vert \right)  \right)  }\int_{0\leq x\leq
av}\left(  \theta_{1}^{2}x^{2}+\theta_{2}^{2}v^{2}\right)  F\left(
\frac{\mathrm{d}x}{v}\right)  \Lambda(\mathrm{d}v),
\]
and when $\theta_{1}\theta_{2}<0$ we have $\theta_{1}\theta_{2}\int_{-va\leq
x\leq0}xF\left(  \frac{\mathrm{d}x}{v}\right)  \geq0$, which gives%
\[
\int_{0\leq v\leq1/\left(  2a\left(  \left\vert \theta_{1}\right\vert
\vee\left\vert \theta_{2}\right\vert \right)  \right)  }\int_{\left\vert
x\right\vert \leq va}\left(  \cos(\theta_{1}x+\theta_{2}v)-1\right)  F\left(
\frac{\mathrm{d}x}{v}\right)  \Lambda(\mathrm{d}v)
\]%
\[
\leq-c\int_{0\leq v\leq1/\left(  2a\left(  \left\vert \theta_{1}\right\vert
\vee\left\vert \theta_{2}\right\vert \right)  \right)  }\int_{-av\leq x\leq
0}\left(  \theta_{1}^{2}x^{2}+\theta_{2}^{2}v^{2}\right)  F\left(
\frac{\mathrm{d}x}{v}\right)  \Lambda(\mathrm{d}v).
\]
Notice that
\[
c\int_{0\leq v\leq1/\left(  2a\left(  \left\vert \theta_{1}\right\vert
\vee\left\vert \theta_{2}\right\vert \right)  \right)  }\int_{0\leq x\leq
av}\theta_{2}^{2}v^{2}F\left(  \frac{\mathrm{d}x}{v}\right)  \Lambda
(\mathrm{d}v)
\]%
\[
=c\left(  F\left(  a\right)  -F\left(  0\right)  \right)  \int_{0\leq
v\leq1/\left(  2a\left(  \left\vert \theta_{1}\right\vert \vee\left\vert
\theta_{2}\right\vert \right)  \right)  }\theta_{2}^{2}v^{2}\Lambda
(\mathrm{d}v).
\]
We get by arguing as on the top of page 968 in Pruitt \cite{pruitt} or in the
remark after the proof of Proposition 6.1 in Buchmann, Maller and Mason
\cite{BMM}, that for some $c_{1}>0$ and $0<\kappa<1$, that whenever $2a\left(
\left\vert \theta_{1}\right\vert \vee\left\vert \theta_{2}\right\vert \right)
\geq1$%
\[
-c\left(  F\left(  a\right)  -F\left(  0\right)  \right)  \theta_{2}^{2}%
\int_{0\leq v\leq1/\left(  2a\left(  \left\vert \theta_{1}\right\vert
\vee\left\vert \theta_{2}\right\vert \right)  \right)  }v^{2}\Lambda
(\mathrm{d}v)\leq-\frac{c_{1}\theta_{2}^{2}}{4a^{2}\left(  \left\vert
\theta_{1}\right\vert \vee\left\vert \theta_{2}\right\vert \right)  ^{2}%
}\left(  2a\left(  \left\vert \theta_{1}\right\vert \vee\left\vert \theta
_{2}\right\vert \right)  \right)  ^{\kappa}.
\]
Next,
\[
-c\int_{0\leq v\leq1/\left(  2a\left(  \left\vert \theta_{1}\right\vert
\vee\left\vert \theta_{2}\right\vert \right)  \right)  }\int_{0\leq x\leq
av}\theta_{1}^{2}x^{2}F\left(  \frac{\mathrm{d}x}{v}\right)  \Lambda
(\mathrm{d}v)
\]%
\[
=-c\theta_{1}^{2}\int_{0\leq x\leq a}u^{2}F\left(  \mathrm{d}u\right)
\int_{0\leq v\leq1/\left(  2a\left(  \left\vert \theta_{1}\right\vert
\vee\left\vert \theta_{2}\right\vert \right)  \right)  }v^{2}\Lambda
(\mathrm{d}v),
\]
which by the previous argument is for some $c_{2}>0$, for $2a\left(
\left\vert \theta_{1}\right\vert \vee\left\vert \theta_{2}\right\vert \right)
\geq1$%
\[
\leq-\frac{c_{2}\theta_{1}^{2}}{\left(  2a\left(  \left\vert \theta
_{1}\right\vert \vee\left\vert \theta_{2}\right\vert \right)  \right)  ^{2}%
}\left(  2a\left(  \left\vert \theta_{1}\right\vert \vee\left\vert \theta
_{2}\right\vert \right)  \right)  ^{\kappa}.
\]
Thus with $c_{3}=c_{1}\wedge c_{2}$,
\[
-c\int_{0\leq v\leq1/\left(  2a\left(  \left\vert \theta_{1}\right\vert
\vee\left\vert \theta_{2}\right\vert \right)  \right)  }\int_{0\leq x\leq
av}(\theta_{1}^{2}x^{2}+\theta_{2}^{2}v^{2})F\left(  \frac{\mathrm{d}x}%
{v}\right)  \Lambda(\mathrm{d}v)
\]%
\[
\leq-c_{3}\left(  \frac{\theta_{1}^{2}}{4a^{2}\left(  \left\vert \theta
_{1}\right\vert \vee\left\vert \theta_{2}\right\vert \right)  ^{2}}%
+\frac{\theta_{2}^{2}}{4a^{2}\left(  \left\vert \theta_{1}\right\vert
\vee\left\vert \theta_{2}\right\vert \right)  ^{2}}\right)  \left(  2a\left(
\left\vert \theta_{1}\right\vert \vee\left\vert \theta_{2}\right\vert \right)
\right)  ^{\kappa}.
\]
Notice that%
\[
\frac{\theta_{1}^{2}}{4a^{2}\left(  \left\vert \theta_{1}\right\vert
\vee\left\vert \theta_{2}\right\vert \right)  ^{2}}+\frac{\theta_{2}^{2}%
}{4a^{2}\left(  \left\vert \theta_{1}\right\vert \vee\left\vert \theta
_{2}\right\vert \right)  ^{2}}\geq\frac{1}{4a^{2}}.
\]
Hence when $\theta_{1}\theta_{2}>0$ and $2a\left(  \left\vert \theta
_{1}\right\vert \vee\left\vert \theta_{2}\right\vert \right)  \geq1$ for some
$c_{4}>0$,
\begin{equation}
-c\int_{0\leq v\leq1/\left(  2a\left(  \left\vert \theta_{1}\right\vert
\vee\left\vert \theta_{2}\right\vert \right)  \right)  }\int_{0\leq x\leq
av}(\theta_{1}^{2}x^{2}+\theta_{2}^{2}v^{2})F\left(  \frac{\mathrm{d}x}%
{v}\right)  \Lambda(\mathrm{d}v)\leq-c_{4}\left(  \left\vert \theta
_{1}\right\vert \vee\left\vert \theta_{2}\right\vert \right)  ^{\kappa
}.\label{e1}%
\end{equation}
The analogous inequality holds when $\theta_{1}\theta_{2}\leq0$ and $2a\left(
\left\vert \theta_{1}\right\vert \vee\left\vert \theta_{2}\right\vert \right)
\geq1$, namely for some $c_{5}>0,$%
\[
\int_{0\leq v\leq1/\left(  2a\left(  \left\vert \theta_{1}\right\vert
\vee\left\vert \theta_{2}\right\vert \right)  \right)  }\int_{\left\vert
x\right\vert \leq va}\left(  \cos(\theta_{1}x+\theta_{2}v)-1\right)  F\left(
\frac{\mathrm{d}x}{v}\right)  \Lambda(\mathrm{d}v)
\]
\begin{equation}
\leq-c\int_{0\leq v\leq1/\left(  2a\left(  \left\vert \theta_{1}\right\vert
\vee\left\vert \theta_{2}\right\vert \right)  \right)  }\int_{-av\leq x\leq
0}\left(  \theta_{1}^{2}x^{2}+\theta_{2}^{2}v^{2}\right)  F\left(
\frac{\mathrm{d}x}{v}\right)  \Lambda(\mathrm{d}v)\leq-c_{5}\left(  \left\vert
\theta_{1}\right\vert \vee\left\vert \theta_{2}\right\vert \right)  ^{\kappa
}.\label{e2}%
\end{equation}
Note that since $0<\kappa<1$ the function $\rho\left(  u\right)  =\left\vert
u\right\vert ^{\kappa}$ is concave on $\left(  0,\infty\right)  $, and thus
\[
\left(  \left\vert \theta_{1}\right\vert \vee\left\vert \theta_{2}\right\vert
\right)  ^{\kappa}\geq\left\vert \frac{\left\vert \theta_{1}\right\vert
+\left\vert \theta_{2}\right\vert }{2}\right\vert ^{\kappa}\geq\frac
{\left\vert \theta_{1}\right\vert ^{\kappa}+\left\vert \theta_{2}\right\vert
^{\kappa}}{2},
\]
which, in combination with (\ref{e1}) and (\ref{e2}), gives for some $d>0,$
whenever $2a\left(  \left\vert \theta_{1}\right\vert \vee\left\vert \theta
_{2}\right\vert \right)  \geq1,$%

\[
\int_{0\leq v\leq1/\left(  2a\left(  \left\vert \theta_{1}\right\vert
\vee\left\vert \theta_{2}\right\vert \right)  \right)  }\int_{\left\vert
x\right\vert \leq va}\left(  \cos(\theta_{1}x+\theta_{2}v)-1\right)  F\left(
\frac{\mathrm{d}x}{v}\right)  \Lambda(\mathrm{d}v)\leq-d\left(  \left\vert
\theta_{1}\right\vert ^{\kappa}+\left\vert \theta_{2}\right\vert ^{\kappa
}\right)  .
\]
\hfill$\Box$ \medskip

The lemma implies that whenever $2a\left(  \left\vert \theta_{1}\right\vert
\vee\left\vert \theta_{2}\right\vert \right)  \geq1$, then for some $d>0$ and
$0<\kappa<1,$%
\[
\left\vert Ee^{\mathrm{i}(\theta_{1}U+\theta_{2}V)}\right\vert \leq\exp\left(
-d\left(  \left\vert \theta_{1}\right\vert ^{\kappa}+\left\vert \theta
_{2}\right\vert ^{\kappa}\right)  \right)  .
\]
As in \cite{pruitt} this allows us to apply the inversion formula for
densities and shows that it may be repeatedly differentiated, from which we
readily infer that $\left(  U,V\right)  $ has a $C^{\infty}$ density when
$EX=0.$ If $EX=\mu\neq0$, the same argument applied to $(U^{\prime},V)=(U-\mu
V,V)$ shows that $(U^{\prime},V)$ has a $C^{\infty}$ density, which by a
simple transformation implies that $\left(  U,V\right)  $ does too.
\hfill$\Box$ \medskip

\textbf{Proof of Corollary \ref{cor1}.} Note that each $V_{t_{k}}/B(t_{k})$ is
an infinitely divisible random variable without a normal component with
L\'{e}vy measure concentrated on $\left(  0,\infty\right)  $ given by
$t_{k}\Lambda\left(  \cdot B(t_{k})\right)  $ with characteristic function
\[
\Psi_{k}\left(  \theta\right)  =\exp\left\{  \mathrm{i}\theta b_{k}+\int
_{0}^{\infty}\left(  e^{\mathrm{i}\theta x}-1-\mathrm{i}\theta x\mathbf{1}%
_{\{0<x\leq1\}}\right)  t_{k}\Lambda(B(t_{k})\mathrm{d}x)\right\}  ,
\]
where
\[
b_{k}=\int_{0}^{1}xt_{k}\Lambda(B(t_{k})\mathrm{d}x).
\]
Since $V_{t_{k}}/B(t_{k})\overset{\mathrm{D}}{\rightarrow}W_{2},$ by
Proposition 7.8 of Sato \cite{sato}, $W_{2}$ is infinitely divisible. Since
$W_{2}$ is necessarily non-negative, it does not have a normal component and
has a L\'{e}vy measure $\Lambda_{0}$ concentrated on $\left(  0,\infty\right)
$. Now by Theorem \ref{2-dim-conv} and its proof, necessarily $\int_{0}%
^{1}x\Lambda_{0}(\mathrm{d}x)<\infty$ and $W_{2}$ has characteristic function
\[
\Psi_{0}\left(  \theta\right)  =\exp\left\{  \mathrm{i}\theta b+\int
_{0}^{\infty}\left(  e^{\mathrm{i}\theta x}-1\right)  \Lambda_{0}%
(\mathrm{d}x)\right\}  ,
\]
where $b\geq0$. By (\ref{Levy-conv-V}) and (\ref{1stmoment-V}) in the proof of
Theorem \ref{2-dim-conv} for any continuity point $v>0$ of $\overline{\Lambda
}_{0}$,
\begin{equation}
t_{k}\overline{\Lambda}(vB(t_{k}))\rightarrow\overline{\Lambda}_{0}%
(v),\quad\text{as }k\rightarrow\infty,\label{v1}%
\end{equation}
and
\begin{equation}
\int_{0}^{v}xt_{k}\Lambda(B(t_{k})\mathrm{d}x)\rightarrow\int_{0}^{v}%
x\Lambda_{0}(\mathrm{d}x)+b,\quad\text{as }k\rightarrow\infty.\label{v2}%
\end{equation}
From (\ref{v2}) we easily get that for any continuity point $v>0$ of
$\overline{\Lambda}_{0}$,%
\begin{equation}
\int_{0}^{v}x^{2}t_{k}\Lambda(B(t_{k})\mathrm{d}x)\rightarrow\int_{0}^{v}%
x^{2}\Lambda_{0}(\mathrm{d}x)=V_{0,2}\left(  v\right)  ,\quad\text{as
}k\rightarrow\infty.\label{v3}%
\end{equation}
(Recall the notation (\ref{V2}).) Now, since $V_{t}$ is in the centered Feller
class, (\ref{cfc}) implies that for some $K>0$
\begin{equation}
\limsup_{k\rightarrow\infty}\frac{v^{2}B^{2}(t_{k})\overline{\Lambda}\left(
vB(t_{k})\right)  }{V_{2}\left(  vB(t_{k})\right)  }\leq K\text{.}\label{K}%
\end{equation}
Note
\[
\frac{v^{2}B^{2}(t_{k})\overline{\Lambda}\left(  vB(t_{k})\right)  }%
{V_{2}\left(  vB(t_{k})\right)  }=\frac{v^{2}t_{k}\overline{\Lambda}%
(vB(t_{k}))}{\int_{0}^{v}x^{2}t_{k}\Lambda(B(t_{k})\mathrm{d}x)},
\]
which by (\ref{v1}) and (\ref{v3}) converges to $v^{2}\overline{\Lambda}%
_{0}(v)/V_{0,2}\left(  v\right)  $ for each continuity point $v>0$ of
$\overline{\Lambda}_{0}$. This with (\ref{K}) implies that
\[
\sup_{v>0}\frac{v^{2}\overline{\Lambda}_{0}(v)}{\int_{0}^{v}x^{2}\Lambda
_{0}(\mathrm{d}x)}\leq K,
\]
so Proposition \ref{prop-fc0} applies. \hfill$\Box$ \medskip

\textbf{Proof of Proposition \ref{sums-th4}.} The proof is a simple adaptation
of the proof of Theorem 4 in \cite{KM}, so we only sketch it here. Putting
\[
B_{t}(k)=\left\{  \frac{\left\vert \sum_{i=2}^{\infty}X_{i}\varphi
(S_{i}/t)\right\vert }{\sum_{i=1}^{\infty}\varphi(S_{i}/t)}\leq\frac
{E|X|}{\sqrt{k}}\right\}  ,
\]
and recalling definition (\ref{AAa}), the conditional version of Chebyshev's
inequality implies that $P\{B_{t}(k)|A_{t}(k^{-1})\}\geq1-1/\sqrt{k}$.
Noticing that on the set $B_{t}(k)\cap A_{t}(k^{-1})$
\[
\Delta_{t}\leq\frac{|X_{1}|}{k}+\frac{E|X|}{\sqrt{k}},
\]
a tightness argument finishes the proof. \hfill$\Box$ \medskip

Now we are ready to prove Theorem \ref{subseq}.

Choose any cdf $F$ in the class $\mathcal{X}$. Corollary \ref{cor1} says
whenever $V_{t}$ is in the centered Feller class at 0 $(\infty)$ then every
subsequential limit law of $U_{t}/V_{t}$, as $t\searrow0$, (as $t\rightarrow
\infty$) has a Lebesgue density on $\mathbb{R}$ and hence is continuous.

Suppose $V_{t}$ is in the Feller class at 0 $(\infty),$ but not in the
centered Feller class at 0 $(\infty)$. In this case Corollary \ref{cor2}
implies that one of the subsequential limits is the constant $EX$ and thus not
continuous. Next Proposition 5.5 in \cite{MM3} in the case $t\searrow0$ and
Proposition 3.2 in \cite{MM4} in the case $t\rightarrow\infty$ show that
whenever $V_{t}$ is not in the Feller class at 0 ($\infty$), that is
\[
\limsup_{t\searrow0\text{ }\left(  t\rightarrow\infty\right)  }\frac
{t^{2}\overline{\Lambda}(t)}{\int_{0}^{t}y^{2}\Lambda(\mathrm{d}y)}=\infty,
\]
and (\ref{sc-condition}) holds, then there exist a subsequence $t_{k}%
\searrow0$ ($t_{k}\rightarrow\infty$), such that (\ref{A-assump}) holds, which
by Corollary \ref{cor3} for any $X$ such that $P\{X=x_{0}\}>0$ for some
$x_{0}$, there exists a subsequence $t_{k}\searrow0$ $(t_{k}\rightarrow
\infty)$ such that $U_{t_{k}}/V_{t_{k}}\overset{\mathrm{D}}{\longrightarrow}%
T$, with $P\{T=x_{0}\}>0$, that is, $T$ is not continuous. Finally, assume
that (\ref{sc-condition}) does not hold, then by Proposition \ref{non-feller}
there exists a subsequence $t_{k}\searrow0$ or $t_{k}\rightarrow\infty$ such
that $U_{t_{k}}/V_{t_{k}}\overset{\mathrm{D}}{\longrightarrow}T,$ with
$P\{T=EX\}>0$, and again $T$ is not continuous. This completes the proof of
Theorem \ref{subseq}.

\section{Appendix}

To finish the proofs of Proposition \ref{prop-rt} and thus Theorem \ref{Th2}
we shall require the following technical result.

\begin{proposition}
\label{P1} Assume that
\begin{equation}
\liminf_{s\searrow0} \frac{s \overline\Lambda(s) }{ \int_{0}^{s}
\overline\Lambda( x ) \mathrm{d}x } >0 \text{,}\label{inf}%
\end{equation}
then
\begin{equation}
ER_{t}=\int_{0}^{\infty}f_{\left(  t\right)  }\left(  u\right)  \mathrm{d}%
u=-t\int_{0}^{\infty}u\Phi^{\prime\prime}\left(  u\right)  e^{-t\Phi\left(
u\right)  }\mathrm{d}u.\label{int}%
\end{equation}

\end{proposition}

\noindent\textbf{Proof.} Clearly for each $u>0$, $f_{T,t}\left(  u\right)
\rightarrow f_{\left(  t\right)  }\left(  u\right)  $, as $T\rightarrow\infty
$. Therefore by Fatou's lemma
\begin{equation}
\int_{0}^{\infty}f_{\left(  t\right)  }\left(  u\right)  \mathrm{d}%
u\leq\liminf_{T\rightarrow\infty}\int_{0}^{\infty}f_{T,t}\left(  u\right)
\mathrm{d}u=\liminf_{T\rightarrow\infty}ER_{t}\left(  T\right)  \leq
1.\label{Tone}%
\end{equation}
Keeping in mind (\ref{conv}) and (\ref{ERT}), this implies that
\[
\int_{0}^{\infty}f_{\left(  t\right)  }\left(  u\right)  \mathrm{d}u\leq
ER_{t}\leq1.
\]
Therefore on account of (\ref{conv}) to prove (\ref{int}) it suffices to
establish (\ref{PL}), as $T\rightarrow\infty.$ One can readily check using
(\ref{delta}) that for some constants $C_{1}>0$ and $C_{2}>0$ and all $u>0$
\[
0\leq-tu\Phi^{\prime\prime}(u)\leq t\left(  C_{1}+u^{-1}C_{2}\right)  .
\]
To see this note that for each $u>0$%
\begin{align*}
-u\Phi^{\prime\prime}(u) &  =u\int_{0}^{\infty}x^{2}e^{-ux}\Lambda
(\mathrm{d}x)\\
&  =\int_{0}^{1}x^{2}ue^{-ux}\Lambda(\mathrm{d}x)+u^{-1}\int_{1}^{\infty}%
u^{2}x^{2}e^{-ux}\Lambda(\mathrm{d}x),\\
&  \leq\max_{0\leq y}ye^{-y}\int_{0}^{1}x\Lambda(\mathrm{d}x)+u^{-1}%
\overline{\Lambda}\left(  1\right)  \max_{0\leq y}y^{2}e^{-y}=:C_{1}%
+u^{-1}C_{2}\text{.}%
\end{align*}
Thus since
\[
f_{T,t}\left(  u\right)  \leq-ut\Phi_{T}^{\prime\prime}(u)\leq-ut\Phi
^{\prime\prime}(u),
\]
we get by the bounded convergence theorem that for all $D>\delta>0$%
\[
\lim_{T\rightarrow\infty}\int_{\delta}^{D}f_{T,t}\left(  u\right)
\mathrm{d}u=\int_{\delta}^{D}f_{\left(  t\right)  }\left(  u\right)
\mathrm{d}u.
\]
Notice that since
\[
\int_{0}^{\infty}f_{\left(  t\right)  }\left(  u\right)  \mathrm{d}u\leq1,
\]
we have%
\[
\lim_{\delta\rightarrow0}\int_{0}^{\delta}f_{\left(  t\right)  }\left(
u\right)  \mathrm{d}u=0\text{ and }\lim_{D\rightarrow\infty}\int_{D}^{\infty
}f_{\left(  t\right)  }(u)\mathrm{d}u=0.
\]
We see now that to complete the verification of (\ref{PL}) we have to show
that
\begin{equation}
\lim_{\delta\rightarrow0}\limsup_{T\rightarrow\infty}\int_{0}^{\delta}%
f_{T,t}\left(  u\right)  \mathrm{d}u=0\label{delta0}%
\end{equation}
and
\begin{equation}
\lim_{D\rightarrow\infty}\limsup_{T\rightarrow\infty}\int_{D}^{\infty}%
f_{T,t}(u)\mathrm{d}u=0.\label{Dinfty}%
\end{equation}
The first condition (\ref{delta0}) is easy to show. Recalling (\ref{fT}),
notice that
\[
f_{T,t}(u)\leq tu\int_{0}^{\infty}\varphi^{2}(s)e^{-u\varphi(s)}\mathrm{d}s,
\]
and so by Fubini
\begin{align*}
\int_{0}^{\delta}f_{T,t}(u)\mathrm{d}u &  \leq t\int_{0}^{\infty}\varphi
^{2}(s)\mathrm{d}s\int_{0}^{\delta}ue^{-u\varphi(s)}\mathrm{d}u\\
&  =t\int_{0}^{\infty}\left[  -\varphi(s)\delta e^{-\delta\varphi
(s)}+(1-e^{-\delta\varphi(s)})\right]  \mathrm{d}s\\
&  =t\left(  \Phi(\delta)-\delta\Phi^{\prime}(\delta)\right)  \leq
t\Phi(\delta),
\end{align*}
which goes to $0$ as $\delta\rightarrow0$ and thus (\ref{delta0})
holds.\smallskip

For the second condition (\ref{Dinfty}), choose $D>0.$ We see that for all
large enough $T>0$%
\begin{equation}
\int_{D}^{\infty}f_{T,t}(u)\mathrm{d}u=\int_{D}^{1/\varphi(T)}f_{T,t}%
(u)\mathrm{d}u+\int_{1/\varphi(T)}^{\infty}f_{T,t}(u)\mathrm{d}u.\label{qq}%
\end{equation}
Recall that
\begin{equation}
f_{T,t}(u)=tu\int_{0}^{T}\varphi^{2}(s)e^{-u\varphi(s)}\mathrm{d}s\exp\left\{
-t\int_{0}^{T}\left(  1-e^{-u\varphi(s)}\right)  \mathrm{d}s\right\}
.\label{f2}%
\end{equation}
We shall first bound the second integral on the right side of (\ref{qq}). For
$u\varphi(T)\geq1$ and keeping mind that $\varphi(s)\geq\varphi(T)$ for
$0<s\leq T$, we have
\[
\exp\left\{  -t\int_{0}^{T}\left(  1-e^{-u\varphi(s)}\right)  \mathrm{d}%
s\right\}  \leq e^{-t(1-e^{-1})T}%
\]
and so
\[
\int_{1/\varphi(T)}^{\infty}f_{T,t}(u)\mathrm{d}u\leq te^{-t(1-e^{-1})T}%
\int_{1/\varphi(T)}^{\infty}u\int_{0}^{T}\varphi^{2}(s)e^{-u\varphi
(s)}\mathrm{d}s\mathrm{d}u.
\]
Using Fubini's theorem the last integral is easy to calculate. We get
\begin{align*}
\int_{1/\varphi(T)}^{\infty}u\int_{0}^{T}\varphi^{2}(s)e^{-u\varphi
(s)}\mathrm{d}s\mathrm{d}u  &  =\int_{0}^{T}\varphi^{2}(s)\mathrm{d}%
s\int_{1/\varphi(T)}^{\infty}ue^{-u\varphi(s)}\mathrm{d}u\\
&  =\int_{0}^{T}\left(  e^{-\varphi(s)/\varphi(T)}+\frac{\varphi(s)}%
{\varphi(T)}e^{-\varphi(s)/\varphi(T)}\right)  \mathrm{d}s\\
&  \leq T\left(  1+\max_{y\geq0}ye^{-y}\right)  \leq2T.
\end{align*}
So we obtain
\begin{equation}
\int_{1/\varphi(T)}^{\infty}f_{T,t}(u)\mathrm{d}u\leq2Tte^{-t(1-e^{-1}%
)T},\label{D1}%
\end{equation}
which tends to $0$ as $T\rightarrow\infty$.\smallskip

\noindent Therefore to complete the verification that (\ref{Dinfty}) holds and
thus (\ref{PL}) we must prove that
\begin{equation}
\lim_{D\rightarrow\infty}\limsup_{T\rightarrow\infty}\int_{D}^{1/\varphi
(T)}f_{T,t}(u)\mathrm{d}u=0.\label{D}%
\end{equation}
We shall bound $f_{T,t}\left(  u\right)  $ in the integral (\ref{D}). Since
$1/u\geq\varphi(T)$, and thus $\overline{\Lambda}\left(  1/u\right)
\leq\overline{\Lambda}\left(  \varphi(T)\right)  \leq T$, we get that the
second factor of $f_{T,t}\left(  u\right)  $ given in (\ref{f2}) is
\[%
\begin{split}
\exp\left\{  -t\int_{0}^{T}\left(  1-e^{-u\varphi(s)}\right)  \mathrm{d}%
s\right\}   &  \leq\exp\left\{  -t\int_{0}^{\overline{\Lambda}\left(
1/u\right)  }\left(  1-e^{-u\varphi(s)}\right)  \mathrm{d}s\right\} \\
&  \leq e^{-t(1-e^{-1})\overline{\Lambda}(1/u)}.
\end{split}
\]
While for the first factor in $f_{T,t}(u)$ given in (\ref{f2}) we use the
simple bound
\[
tu\int_{0}^{T}\varphi^{2}(s)e^{-u\varphi(s)}\mathrm{d}s\leq tu\int_{0}%
^{\infty}\varphi^{2}(s)e^{-u\varphi(s)}\mathrm{d}s=:t\psi_{\Lambda}\left(
u\right)  .
\]
We see that%
\[%
\begin{split}
\int_{D}^{1/\varphi(T)}f_{T,t}(u)\mathrm{d}u  &  \leq t\int_{D}^{1/\varphi
(T)}\psi_{\Lambda}\left(  u\right)  e^{-t(1-e^{-1})\overline{\Lambda}%
(1/u)}\mathrm{d}u\\
&  \leq t\int_{D}^{\infty}\psi_{\Lambda}\left(  u\right)  e^{-t(1-e^{-1}%
)\overline{\Lambda}(1/u)}\mathrm{d}u.
\end{split}
\]
Clearly (\ref{Dinfty}) holds whenever for all $\gamma>0$,
\begin{equation}
\int_{1}^{\infty}\psi_{\Lambda}\left(  u\right)  e^{-\gamma\overline{\Lambda
}(1/u)}\mathrm{d}u<\infty.\label{levy-cond}%
\end{equation}

\begin{lemma}
\label{L1} Whenever (\ref{inf}) is satisfied, then for all $\gamma>0$,
(\ref{levy-cond}) holds.
\end{lemma}

\noindent\textbf{Proof.} Recall the definition (\ref{def-f(u)}). Notice that
by (\ref{Tone}) for all $t>0$%
\begin{equation}
\int_{0}^{\infty}f_{\left(  t\right)  }\left(  u\right)  \mathrm{d}%
u<\infty.\label{t}%
\end{equation}
Write%
\[
\int_{0}^{\infty}(1-e^{-u\varphi(s)})\mathrm{d}s=\int_{0}^{1/u}(1-e^{-ux}%
)\Lambda(\mathrm{d}x)+\int_{1/u}^{\infty}(1-e^{-ux})\Lambda\left(
\mathrm{d}x\right)  .
\]
We see that%
\[
\int_{1/u}^{\infty}(1-e^{-ux})\Lambda\left(  \mathrm{d}x\right)  \leq
\overline{\Lambda}\left(  1/u\right)
\]
and%
\[%
\begin{split}
\int_{0}^{1/u}(1-e^{-ux})\Lambda\left(  \mathrm{d}x\right)   &  =-\left(
1-e^{-1}\right)  \overline{\Lambda}(1/u)+\int_{0}^{1/u}u\overline{\Lambda
}\left(  x\right)  e^{-ux}\mathrm{d}x\\
&  \leq\int_{0}^{1/u}u\overline{\Lambda}\left(  x\right)  e^{-ux}%
\mathrm{d}x\leq u\int_{0}^{1/u}\overline{\Lambda}\left(  x\right)
\mathrm{d}x.
\end{split}
\]
By assumption (\ref{inf}) for some $\eta>0$ for all $u$ large%
\begin{equation}
u\int_{0}^{1/u}\overline{\Lambda}\left(  x\right)  \mathrm{d}x\leq
\eta\overline{\Lambda}(1/u).\label{tail}%
\end{equation}
This implies that%
\[
t\int_{0}^{\infty}(1-e^{-u\varphi(s)})\mathrm{d}s\leq\left(  1+\eta\right)
t\overline{\Lambda}(1/u).
\]
Thus for all large enough $D>0$ and all $t>0$
\[
\int_{D}^{\infty}f_{\left(  t\right)  }\left(  u\right)  \mathrm{d}u\geq
\int_{D}^{\infty}t\psi_{\Lambda}\left(  u\right)  \exp\left\{  -\left(
1+\eta\right)  t\overline{\Lambda}(1/u)\right\}  \mathrm{d}u,
\]
and hence since (\ref{t}) holds for all $t>0$, we get that for all $\gamma>0
$, (\ref{levy-cond}) is satisfied. \hfill$\Box\medskip$

\noindent We see from Lemma \ref{L1} that (\ref{levy-cond}) holds whenever
assumption (\ref{inf}) is satisfied and thus by the arguments preceding the
lemma the limit (\ref{PL}) is valid.\textbf{\ }This completes the proof of
Proposition \ref{P1}. \hfill$\Box$

\subsection{Return to the proofs of Proposition \ref{prop-rt} and Theorem
\ref{Th2}}

We shall now finish the proof of Proposition \ref{prop-rt}. To do this we
shall need three more lemmas. Let $X_{t}$ be a subordinator with canonical
measure $\Lambda$. Assume that $X_{t}$ is without drift. Define
\[
I(x)=\int_{0}^{x}\overline{\Lambda}(y)\mathrm{d}y.
\]
We give a criterion for subsequential relative stability of $X$ at $0$.

\begin{lemma}
\label{th1} Let $X$ be a driftless subordinator with $\overline{\Lambda
}(0+)>0$. There are nonstochastic sequences $t_{k}\downarrow0$ and $B_{k}>0$,
such that, as $k\rightarrow\infty$,
\begin{equation}
\frac{X(t_{k})}{B_{k}}\overset{\mathrm{P}}{\longrightarrow}1\label{0.1}%
\end{equation}
if and only if
\begin{equation}
\liminf_{x\downarrow0}\frac{x\overline{\Lambda}(x)}{I(x)}=0.\label{0.2}%
\end{equation}

\end{lemma}

\noindent\textbf{Proof.} From the convergence criteria for subordinators, e.g.
part (ii) of Theorem 15.14 of \cite{kallenberg}, p. 295, (\ref{0.1}%
)\textbf{\ }is equivalent to
\begin{equation}
\lim_{t_{k}\rightarrow0}t_{k}\overline{\Lambda}(xB_{k})=0\text{ for every
}x>0\text{ }\mathrm{and}\text{ }\lim_{t_{k}\rightarrow0}t_{k}\int_{0}%
^{1}x\Lambda(\mathrm{d}B_{k}x)=1.\label{kall}%
\end{equation}
Noting that $I(x)=\int_{0}^{x}y\Lambda(\mathrm{d}y)+x\overline{\Lambda}(x)$,
we see that (\ref{kall})\textbf{\ }implies
\begin{equation}
t_{k}B_{k}^{-1}I(B_{k})=t_{k}B_{k}^{-1}\int_{0}^{B_{k}}x\Lambda(\mathrm{d}%
x)+t_{k}\overline{\Lambda}(B_{k})\rightarrow1\text{,}\label{0.3}%
\end{equation}
and clearly (\ref{0.3}) and $t_{k}\overline{\Lambda}(B_{k})\rightarrow0$ imply
(\ref{0.2}). (Note that necessarily $B_{k}\rightarrow0$.)\medskip\ 

\noindent Conversely, let (\ref{0.2}) hold and choose a subsequence
$c_{k}\rightarrow0$ as $k\rightarrow\infty$ such that
\[
\lim_{k\rightarrow\infty}\frac{c_{k}\overline{\Lambda}(c_{k})}{I(c_{k})}=0.
\]
Define
\[
t_{k}:=\sqrt{\frac{c_{k}}{\overline{\Lambda}(c_{k})I(c_{k})}}.
\]
Then
\[
\lim_{k\rightarrow\infty}t_{k}\overline{\Lambda}(c_{k})=\lim_{k\rightarrow
\infty}\sqrt{\frac{c_{k}\overline{\Lambda}(c_{k})}{I(c_{k})}}=0,
\]
and
\[
\lim_{k\rightarrow\infty}\frac{t_{k}I(c_{k})}{c_{k}}=\lim_{k\rightarrow\infty
}\sqrt{\frac{I(c_{k})}{c_{k}\overline{\Lambda}(c_{k})}}=\infty.
\]
Then set $B_{k}:=t_{k}I(c_{k})$, so $\lim_{k\rightarrow\infty}B_{k}=0$ and
$\lim_{k\rightarrow\infty}B_{k}/c_{k}=\infty$. Given $x>0$ choose $k$ so large
that $xB_{k}\geq c_{k}$. Then
\begin{equation}
t_{k}\overline{\Lambda}(xB_{k})\leq t_{k}\overline{\Lambda}(c_{k}%
)\rightarrow0.\label{0.4}%
\end{equation}
Furthermore, by writing
\[
\frac{t_{k}I(B_{k})}{B_{k}}=\frac{t_{k}I(c_{k})}{B_{k}}+\frac{t_{k}\left(
I(B_{k})-I(c_{k})\right)  }{B_{k}}=1+\frac{t_{k}\left(  I(B_{k})-I(c_{k}%
)\right)  }{B_{k}}%
\]
and noting that for all large $k$
\[
0\leq\frac{t_{k}\left(  I(B_{k})-I(c_{k})\right)  }{B_{k}}\leq\frac{B_{k}%
t_{k}\overline{\Lambda}(c_{k})}{B_{k}}\rightarrow0,
\]
we also have $t_{k}B_{k}^{-1}I(B_{k})\rightarrow1$ and thus by (\ref{0.4}) and
the identity in (\ref{0.3})
\[
\lim_{t_{k}\rightarrow0}t_{k}\int_{0}^{1}x\Lambda(\mathrm{d}B_{k}x)=1
\]
which in combination with (\ref{0.4}) implies (\ref{0.1}), by (\ref{kall}).
\hfill$\Box\medskip$

To continue we need the following lemma from \cite{MM}.

\begin{lemma}
\label{psi-lemma} Let $\Psi$ be a non-negative measurable real valued function
defined on $(0,\infty)$ satisfying
\[
\int_{0}^{\infty}\Psi\left(  y\right)  \mathrm{d}y<\infty.
\]
Then
\begin{equation}
E \left(  \sum_{i=1}^{\infty}\Psi\left(  S_{i}\right)  \right)  =\int
_{0}^{\infty}\Psi\left(  y\right)  \mathrm{d}y\label{ee1}%
\end{equation}
and $\lim_{n\rightarrow\infty}E\left(  \sum_{i=n}^{\infty}\Psi\left(
S_{i}\right)  \right)  =0.$
\end{lemma}

\begin{lemma}
\label{lemma-rt-inf} (i) Assume that (\ref{rt}) holds as $t \searrow0$ with
$\beta< 1$. Then (\ref{inf}) holds.

(ii) Assume that (\ref{rt}) holds as $t \to\infty$ with $\beta< 1$. Then
without loss of generality we can assume that (\ref{inf}) holds.
\end{lemma}

\noindent\textbf{Proof.} (i) We shall show that (\ref{rt}) implies
(\ref{inf}). Assume on the contrary that (\ref{inf}) does not hold. Then,
since $V_{t}$ is a driftless subordinator by Lemma \ref{th1} for some
sequences $B_{k}>0$, $t_{k}\downarrow0$, $V_{t_{k}}/B_{k}\overset{\mathrm{P}%
}{\rightarrow}1.$ By Proposition \ref{prop1-repr} the infinite sum $\sum
_{i=1}^{\infty}\varphi\left(  \frac{S_{i}}{t}\right)  $ is equal in
distribution to $V_{t}$ and $\sum_{i=1}^{\infty}\varphi^{2}\left(  \frac
{S_{i}}{t}\right)  $ is equal in distribution to the subordinator $W_{t}$ with
L\'{e}vy measure $\Lambda_{2}$ on $\left(  0,\infty\right)  $ defined by
\[
\overline{\Lambda}_{2}\left(  x\right)  =\overline{\Lambda}\left(  \sqrt
{x}\right)  .
\]
From (\ref{kall}) in the proof of Lemma \ref{th1} above
\begin{equation}
t_{k}\overline{\Lambda}\left(  xB_{k}\right)  \rightarrow0\text{ and }\int
_{0}^{1}t_{k}\overline{\Lambda}\left(  xB_{k}\right)  \mathrm{d}%
x\rightarrow1,\label{twoL}%
\end{equation}
with $t_{k}\rightarrow0$ and $B_{k}\rightarrow0$. Thus we easily see that
\[
t_{k}\overline{\Lambda}_{2}\left(  xB_{k}^{2}\right)  =t_{k}\overline{\Lambda
}\left(  \sqrt{x}B_{k}\right)  \rightarrow0
\]
and
\[
\int_{0}^{1}t_{k}\overline{\Lambda}_{2}\left(  xB_{k}^{2}\right)
\mathrm{d}x=\int_{0}^{1}t_{k}\overline{\Lambda}\left(  \sqrt{x}B_{k}\right)
\mathrm{d}x=2\int_{0}^{1}yt_{k}\overline{\Lambda}\left(  yB_{k}\right)
\mathrm{d}y,
\]
which for any $0<\delta<1$ is
\[
\leq2\delta\int_{0}^{1}t_{k}\overline{\Lambda}\left(  xB_{k}\right)
\mathrm{d}x+2\int_{\delta}^{1}t_{k}\overline{\Lambda}\left(  xB_{k}\right)
\mathrm{d}x.
\]
Clearly by (\ref{twoL})
\[
\limsup_{k\rightarrow\infty}\left(  2\delta\int_{0}^{1}t_{k}\overline{\Lambda
}\left(  xB_{k}\right)  \mathrm{d}x+2\int_{\delta}^{1}t_{k}\overline{\Lambda
}\left(  xB_{k}\right)  \mathrm{d}x\right)  =2\delta.
\]
Thus since $0<\delta<1$ can be made arbitrarily small we get
\[
\lim_{k\rightarrow\infty}\int_{0}^{1}t_{k}\overline{\Lambda}_{2}\left(
xB_{k}^{2}\right)  \mathrm{d}x=0.
\]
Hence applying Theorem 15.14 on page 295 of \cite{kallenberg}, we get
$W_{t_{k}}/B_{k}^{2}\overset{\mathrm{P}}{\rightarrow}0$ and thus
\[
R_{t_{k}}\overset{\mathrm{D}}{=}W_{t_{k}}/\left(  V_{t_{k}}\right)
^{2}\overset{\mathrm{P}}{\rightarrow}0,
\]
which since $R_{t_{k}}\leq1$ implies $ER_{t_{k}}\rightarrow0$, as
$t_{k}\downarrow0$, which clearly contradicts to (\ref{rt}). So we have
(\ref{inf}) in this case. \smallskip

(ii) We shall first assume that
\begin{equation}
\int_{0}^{\infty}\varphi\left(  u\right)  \mathrm{d}u=\infty,\label{infin}%
\end{equation}
which by (\ref{delta}) implies
\begin{equation}
\int_{0}^{1}\varphi\left(  u\right)  \mathrm{d}u=\infty.\label{phh2}%
\end{equation}
Set%
\[
V\left(  t\right)  :=\sum_{i=1}^{\infty}\varphi\left(  \frac{S_{i}}{t}\right)
\mathbf{1}\left\{  \frac{S_{i}}{t}\leq1\right\}  \text{ and }\overline
{V}\left(  t\right)  :=\sum_{i=1}^{\infty}\varphi\left(  \frac{S_{i}}%
{t}\right)  \mathbf{1}\left\{  \frac{S_{i}}{t}>1\right\}  .
\]
We see that%
\[
V\left(  t\right)  \geq\sum_{k=1}^{\infty}\varphi\left(  2^{-k+1}\right)
\sum_{i=1}^{\infty}\mathbf{1}\left\{  2^{-k}<\frac{S_{i}}{t}\leq
2^{-k+1}\right\}  .
\]
Now for each fixed $L\geq1,$ as $t\rightarrow\infty$,
\[
t^{-1}\sum_{k=2}^{L+1}\left(  \varphi\left(  2^{-k+1}\right)  \sum
_{i=1}^{\infty}\mathbf{1}\left\{  2^{-k}<\frac{S_{i}}{t}\leq2^{-k+1}\right\}
\right)  \overset{\mathrm{P}}{\rightarrow}\sum_{k=1}^{L}\varphi\left(
2^{-k}\right)  2^{-k-1}\geq2^{-1}\int_{2^{-L}}^{1}\varphi(u)\mathrm{d}u.
\]
Thus since $L\geq1$ can be made arbitrarily large, on account of
(\ref{phh2}),
\begin{equation}
t^{-1}V\left(  t\right)  \overset{\mathrm{P}}{\rightarrow}\infty,\text{ as
}t\rightarrow\infty\text{.}\label{L11}%
\end{equation}
Next, using (\ref{ee1}), we get
\[
t^{-1}E\overline{V}\left(  t\right)  =t^{-1}\int_{t}^{\infty}\varphi
(y/t)\mathrm{d}y=\int_{1}^{\infty}\varphi(u)\mathrm{d}u<\infty\text{,}%
\]
which implies that
\begin{equation}
t^{-1}\overline{V}\left(  t\right)  =O_{P}\left(  1\right)  ,\text{ as
}t\rightarrow\infty.\label{L22}%
\end{equation}
Hence by (\ref{L11}) and (\ref{L22})
\begin{equation}
\overline{V}\left(  t\right)  /V_{t}\overset{\mathrm{P}}{\rightarrow}0,\text{
as }t\rightarrow\infty.\label{L33}%
\end{equation}
We get then that
\begin{equation}
V_{t}=\sum_{i=1}^{\infty}\varphi\left(  \frac{S_{i}}{t}\right)  =V\left(
t\right)  \left(  1+o\left(  1\right)  \right)  ,\text{ as }t\rightarrow
\infty\text{.}\label{vt}%
\end{equation}
Now set%
\[
W_{t}:=\sum_{i=1}^{\infty}\varphi^{2}\left(  \frac{S_{i}}{t}\right)  \text{,
}W\left(  t\right)  :=\sum_{i=1}^{\infty}\varphi^{2}\left(  \frac{S_{i}}%
{t}\right)  \mathbf{1}\left\{  \frac{S_{i}}{t}\leq1\right\}
\]%
\[
\text{ and }\overline{W}\left(  t\right)  :=\sum_{i=1}^{\infty}\varphi
^{2}\left(  \frac{S_{i}}{t}\right)  \mathbf{1}\left\{  \frac{S_{i}}%
{t}>1\right\}  .
\]
Clearly
\[
t^{-1}E\overline{W}\left(  t\right)  =t^{-1}\int_{t}^{\infty}\varphi
^{2}(y/t)\mathrm{d}y=\int_{1}^{\infty}\varphi^{2}(u)\mathrm{d}u<\infty\text{,}%
\]
which says that $t^{-1}\overline{W}\left(  t\right)  =O_{P}\left(  1\right)  $
as $t\rightarrow\infty$. Hence by (\ref{L33}), $\overline{W}\left(  t\right)
/V_{t}\overset{\mathrm{P}}{\rightarrow}0$ as $t\rightarrow\infty$, which when
combined with (\ref{vt}) gives
\begin{equation}
R_{t}=\frac{W_{t}}{V_{t}^{2}}=\frac{W\left(  t\right)  }{V^{2}\left(
t\right)  }+o_{P}(1)\text{, as }t\rightarrow\infty.\label{rp}%
\end{equation}
Notice that $V(t)$ is a L\'{e}vy process with canonical measure $\Lambda_{1} $
defined via
\[
\overline{\Lambda}_{1}\left(  x\right)  =\overline{\Lambda}\left(  x\right)
,\text{ for }x\geq\varphi\left(  1\right)  \text{, and }\overline{\Lambda}%
_{1}\left(  x\right)  =\overline{\Lambda}\left(  \varphi\left(  1\right)
\right)  \text{ for }0<x<\varphi\left(  1\right)  .
\]
Set $\varphi_{1}(s)=\varphi(s)\mathbf{1}\{s<1\}$. Note that we have
\[
\varphi_{1}\left(  s\right)  =\sup\left\{  y:\overline{\Lambda}_{1}%
(y)>s\right\}  ,\text{ }s>0,
\]
where the supremum of the empty set is taken as 0. Let $R_{t}^{\left(
1\right)  }$ be defined as $R_{t}$ with $\varphi$ replaced by $\varphi_{1}$,
that is,
\[
R_{t}^{\left(  1\right)  }=\frac{W\left(  t\right)  }{\left(  V\left(
t\right)  \right)  ^{2}}=\frac{\sum_{i=1}^{\infty}\varphi_{1}^{2}\left(
\frac{S_{i}}{t}\right)  }{\left(  \sum_{i=1}^{\infty}\varphi_{1}\left(
\frac{S_{i}}{t}\right)  \right)  ^{2}}.
\]
Since $R_{t}\left(  1\right)  =R_{t}^{\left(  1\right)  }$, we see by formula
(\ref{ERT}) that
\begin{equation}
ER_{t}^{\left(  1\right)  }=\int_{0}^{\infty}f_{1,t}\left(  u\right)
\mathrm{d}u.\label{int2}%
\end{equation}
Next from (\ref{rp}), we get $R_{t}^{\left(  1\right)  }-R_{t}\overset
{\mathrm{P}}{\rightarrow}0$, as $t\rightarrow\infty$, which implies that
\[
\lim_{t\rightarrow\infty}ER_{t}=\lim_{t\rightarrow\infty}ER_{t}^{\left(
1\right)  }.
\]
Clearly the tail behavior conclusions about $\Lambda_{1}(x)$, as
$x\rightarrow\infty$, will be identical to those for $\Lambda(x)$, as
$x\rightarrow\infty$. Moreover, since $\overline{\Lambda}_{1}(0+)$ is finite
(\ref{inf}) trivially holds for $\Lambda_{1}$. Therefore in our proof in the
case $t\rightarrow\infty$ we can without loss of generality assume that
(\ref{inf}) is satisfied. \smallskip

\noindent The case $\mu:=$ $\int_{0}^{\infty}\varphi\left(  u\right)
\mathrm{d}u<\infty$ cannot occur when $\beta<1$ in (\ref{rt}). In this case it
is easily checked that
\[
t\overline{\Lambda}\left(  x\mu t\right)  \rightarrow0\text{ for all
}x>0\text{ and }\int_{0}^{1}t\overline{\Lambda}\left(  x\mu t\right)
\mathrm{d}x\rightarrow1.
\]
Therefore by proceeding exactly as above we get that $ER_{t}\rightarrow0$ as
$t\rightarrow\infty$, which forces $\beta=1$. $\hspace*{10pt}$ \hfill$\Box$

Returning to the proof of Proposition \ref{prop-rt}, in the case $t\searrow0
$, Lemma \ref{lemma-rt-inf} shows that the assumption of Proposition \ref{P1}
holds and, in the case $t\rightarrow\infty$, it says that we can assume
without loss of generality that it is satisfied. This completes the proof of
Proposition \ref{prop-rt} and hence that of Theorem \ref{Th2}.$\hspace*{10pt}
$ \hfill$\Box$ \smallskip

\noindent\textbf{Acknowledgement }The authors are grateful to Ross Maller for
Lemma \ref{th1}. Also Jan Rosi\'{n}ski kindly pointed out to us a more
efficient way to prove Proposition \ref{prop1-repr} than our original proof.
The second named author would like to thank the Bolyai Institute of Szeged
University for their hospitality, where this paper was partially written. Our
paper also benefited by useful comments from the referee.

\end{document}